# ASYMPTOTIC RESULTS FOR MAXIMUM LIKELIHOOD ESTIMATORS IN JOINT ANALYSIS OF REPEATED MEASUREMENTS AND SURVIVAL TIME[1]

BY DONGLIN ZENG AND JIANWEN CAI

*University of North Carolina at Chapel Hill*

Maximum likelihood estimation has been extensively used in the joint analysis of repeated measurements and survival time. However, there is a lack of theoretical justification of the asymptotic properties for the maximum likelihood estimators. This paper intends to fill this gap. Specifically, we prove the consistency of the maximum likelihood estimators and derive their asymptotic distributions. The maximum likelihood estimators are shown to be semiparametrically efficient.

**1. Introduction.** Joint analysis of both repeated measurements and survival time has received much attention in the last decade. The motivation of such analysis arises from many medical studies. For example, in an HIV study, the progression of CD4 cell counts in HIV patients and the time to patients' death are of interest. Three different types of questions can be asked. One may wish to know the effect of a particular factor, such as age at the entry, on both the progression of CD4 cell counts and the risk of death. Interest may also arise in studying the longitudinal pattern of CD4 cell counts over a time period; however, the longitudinal path can be truncated due to drop-out or death. In another analysis, one may focus on how the actual CD4 cell count predicts the risk of HIV-related disease, where the true CD4 cell count is often missing or measured with error.

To answer these or similar questions in other medical studies, joint models for repeated measurements and survival time have been proposed. In general, a mixed-effects model (Chapter 9 in [7]) with normal random effects is used to model repeated measurements, while a proportional intensity model [2] is used to model the hazard function of survival time. The covariates in

Received July 2003; revised November 2004.
[1]Supported in part by NIH Grant R01 HL69720.

the models can be subjects' baseline variables, study time or error-free CD4 cell counts, and so on. Random effects are used in both the mixed model and the proportional hazards model to account for the dependence between repeated measurements and survival time due to unobserved heterogeneity. In some of the literature, such a joint model is described as either a selection model or a pattern-mixture model, depending upon how they are derived. When the conditional distribution of survival time given repeated measurements is modeled, the derived joint model is called a selection-model; when the conditional distribution of repeated measurement given survival time is modeled, the derived joint model is called a pattern-mixture model. Selection models have been studied by many authors in different contexts, for example, Tsiatis, DeGruttola and Wulfsohn [23], Wulfsohn and Tsiatis [27] and Xu and Zeger [28, 29]. On the other hand, Wu and Carroll [26], Wu and Bailey [25] and Hogan and Laird [12] proposed pattern-mixture models. In some studies where repeated measurements are considered as an internal covariate predicting risk of survival time, joint analysis is regarded as a missing covariate problem or measurement error problem in a proportional hazards model. For instance, Chen and Little [5] considered the missing covariate problem in the proportional hazards model, although the covariates there were assumed to be time-independent. One referee drew our attention to the paper by Dupuy, Grama and Mesbah [8]. In their paper, the authors stated the asymptotic results for the proportional hazards model with time-dependent covariate, where the covariate was modeled using a parametric distribution instead of a mixed-effects model and the covariate was assumed to be measured without error.

In most of the joint analysis literature, nonparametric maximum likelihood estimation has been proposed (e.g., [23, 27]). Here, nonparametric maximum likelihood estimation means that the nuisance parameter, for example, the cumulative baseline hazard function in the proportional hazards model, can be a function with jumps at some discrete observations (for a complete review of the nonparametric maximum likelihood estimation in a proportional hazards model, refer to [15]). Computationally, the EM algorithm [6] has often been used to calculate the maximum likelihood estimates, where random effects are treated as missing.

However, although the maximum likelihood estimates have been shown to perform well in numerical studies [13], theoretical justification of the asymptotic properties of the maximum likelihood estimates has not been well established, except for Chen and Little [5], who thoroughly studied this issue in a proportional hazards model with missing time-independent covariates. Therefore, in this paper we aim to derive the asymptotic properties of the maximum likelihood estimators in the joint models. Specifically, we rigorously prove the consistency of the maximum likelihood estimators and derive their asymptotic distributions. Our theoretical results further confirm that nonparametric maximum likelihood estimation, which has been proposed in the literature [12, 23, 27], provides efficient estimation. Additionally, we show that the profile likelihood function can be used to give a consistent estimator for the asymptotic variance of the regression coefficients.

The structure of this paper is as follows. A general frame for modeling both repeated measurements and survival time is given in Section 2. The EM algorithm for maximum likelihood estimation is briefly described afterward. Section 3 gives our main results on the asymptotic properties of the maximum likelihood estimators. The proofs for the main theorems are given in Sections 4 and 5. Some technical details are provided in the Appendix.

## 2. Maximum likelihood estimation in joint models.

2.1. *Models and assumptions.* In the joint analysis of repeated measurements and survival times, repeated measurements are considered as the realizations of a certain marker process (usually an internal covariate, Section 6.3.2 in [16]) at finite time points. We use $Y(t)$ to denote the value of such a marker process at time $t$ and we introduce a counting process $R(t)$ (cf. II.1 of [1]) which is right-continuous and only jumps at time $t$ where a measurement is taken. Furthermore, we use $N_T(t)$ and $N_C(t)$ to denote the counting processes generated by survival time $T$ and censoring time $C$, respectively; that is, $N_T(t) = I(T \leq t)$ and $N_C(t) = I(C \leq t)$, where $I(\cdot)$ is the indicator function. Both $Y(t)$ and $T$ are outcome variables of interest.

One essential assumption in all of the joint models proposed in the literature is that the association between the marker process and the survival time is due to observed covariate processes such as baseline information or study time, and so on, which are denoted by $\mathcal{X}(t)$, and unobserved subject-specific effects, which are denoted by $\mathbf{a}$. However, the previous statement is vague, and so we will provide more rigorous assumptions in the following. To do that, we introduce more notation: we denote $\bar{H}_{\mathcal{X}}(t)$ as the longitudinal covariate history prior to time $t$ and denote $\bar{H}_Y(t)$ as the longitudinal response history prior to time $t$; that is, $\bar{H}_{\mathcal{X}}(t) = \{\mathcal{X}(s): s < t\}$ and $\bar{H}_Y(t) = \{Y(s): s < t\}$. Additionally, we denote $\tau$ as the end time of the study. Then we impose the following model assumptions:

(A.1) Random effect $\mathbf{a}$ follows a multivariate normal distribution with mean zero and covariance $\boldsymbol{\Sigma}_a$.
(A.2) For any $t \in [0, \tau]$, the covariate process $\mathcal{X}(t)$ is fully observed and conditional on $\mathbf{a}$, $\bar{H}_{\mathcal{X}}(t), \bar{H}_Y(t)$ and $T \geq t$, the distribution of $\mathcal{X}(t)$ depends only on $\bar{H}_{\mathcal{X}}(t)$. Moreover, with probability one, $\mathcal{X}(t)$ is continuously differentiable in $[0, \tau]$ and $\max_{t \in [0,\tau]} \|\mathcal{X}'(t)\| < \infty$, where $\|\cdot\|$ denotes the Euclidean norm in real space and $\mathcal{X}'(t)$ denotes the derivative of $\mathcal{X}(t)$ with respect to $t$.

(A.3) Conditional on **a**, $\bar{H}_\mathcal{X}(t)$, $\bar{H}_Y(t)$, $\mathcal{X}(t)$ and $T \geq t$, the intensity of the counting process $N_T(t)$ at time $t$ is equal to

$$\lambda(t) \exp\{(\boldsymbol{\phi} \circ \widetilde{\mathbf{W}}(t))^T \mathbf{a} + \mathbf{W}(t)^T \boldsymbol{\gamma}\}, \tag{1}$$

where $\mathbf{W}(t)$ and $\widetilde{\mathbf{W}}(t)$ are sub-processes of $\mathcal{X}(t)$ and $\boldsymbol{\phi}$ is a constant vector of the same dimension as $\widetilde{\mathbf{W}}(t)$. For any vectors $v_1$ and $v_2$ of the same dimension, $v_1 \circ v_2$ is the vector obtained by the component-wise product of $v_1$ and $v_2$.

(A.4) Conditional on **a**, $\bar{H}_\mathcal{X}(t)$, $\bar{H}_Y(t)$, $\mathcal{X}(t)$ and $T > t$, the marker process $Y(t)$ satisfies

$$Y(t) = \mathbf{X}(t)^T \boldsymbol{\beta} + \widetilde{\mathbf{X}}(t)^T \mathbf{a} + \varepsilon(t), \tag{2}$$

where $\mathbf{X}(t)$ and $\widetilde{\mathbf{X}}(t)$ are sub-processes of $\mathcal{X}(t)$ and $\varepsilon(t)$ is a white noise process with variance $\sigma_y^2$.

(A.5) Conditional on **a**, $\bar{H}_\mathcal{X}(T)$, $\bar{H}_Y(T)$, $\mathcal{X}(T)$ and $T$, the intensity of the counting process $N_C(t)$ at time $t$ depends only on $\bar{H}_\mathcal{X}(t)$ and $\mathcal{X}(t)$ for any $t < T$.

(A.6) Conditional on **a**, $\bar{H}_\mathcal{X}(T)$, $\bar{H}_Y(T)$, $\mathcal{X}(T)$, $T$ and $C$, the intensity of the counting process $R(t)$ depends only on $\bar{H}_\mathcal{X}(t)$ and $\mathcal{X}(t)$ for any $t < T \wedge C$.

REMARK 2.1. The structure models (1) and (2) cover most of the joint models proposed in the literature. For example, we consider a clinical trial for HIV-patients with two treatment arms and covariates measured at the entry, such as age and gender and so on, in addition to the recording of CD4 counts along the follow-up. If we are interested in studying the simultaneous effects of the treatments on both CD4 count and survival time, after adjusting for other covariates at the entry and subject-specific effects, we can choose each of $\widetilde{\mathbf{W}}(t)$, $\mathbf{W}(t)$, $\mathbf{X}(t)$ and $\widetilde{\mathbf{X}}(t)$ to include treatment variable and/or covariates measured at the entry. Furthermore, we can use $t$ as a covariate in models (1) and (2) to study time-dependent effects. In another situation, if $Y(t)$ is considered as a marker process subject to measurement error which predicts risk of death (e.g., [23]), then one would use the error-free covariate for $Y(t)$, which equals $\mathbf{X}(t)^T \boldsymbol{\beta} + \widetilde{\mathbf{X}}(t)^T \mathbf{a}$ based on (2), as one predictor with a coefficient $\mu$ in the proportional hazards model (1). Some re-arrangement shows that the obtained hazards model is a special case of the expression (1), where $\widetilde{\mathbf{W}}(t) = \widetilde{\mathbf{X}}(t)$, $\mathbf{W}(t) = \mathbf{X}(t)$, $\boldsymbol{\phi}$ is the vector with each component being $\mu$ and $\boldsymbol{\gamma} = \mu \boldsymbol{\beta}$.

REMARK 2.2. In assumptions (A.5) and (A.6), we implicitly assume that there exist some appropriate measures such that the intensity functions of $N_C(t)$ and $R(t)$ exist. In the missingness context, (A.5) and (A.6) are special cases of the coarsening at random assumption [10]. In fact, we can allow (A.5) and (A.6) to be even more general; for example, the intensities of $N_C(t)$ and $R(t)$ can depend on the observed history of repeated measurements. Since this generality will not change our subsequent arguments, we choose to work with the current simpler assumptions. Assumptions (A.5) and (A.6) can be further simplified under some special situations: when $\mathcal{X}(t)$ is equal to $g(\mathbf{X};t)$, where $g(\cdot)$ is a deterministic function and $\mathbf{X}$ are random variables at time zero, then assumptions (A.5) and (A.6) can be replaced with the following assumptions: conditional on $\mathbf{X}$, $C$ and $R(t)$ are independent of $T$, $\bar{H}_Y(T)$ and $\mathbf{a}$. Such assumptions have been used in [23]. Another special situation when (A.6) is satisfied is when repeated measurements are obtained on fixed schedules, as seen in many cohort studies of clinical experiments.

Under assumptions (A.1)–(A.6), the conditional distribution function for the right-censored event time ($Z = T \wedge C, \Delta = I(T \le C)$), the repeated measurement process $\{(Y(t)I(R(t) - R(t-) > 0), R(t)) : t < Z\}$ and $\{\mathcal{X}(t) : t \le Z\}$, given random effects $\mathbf{a}$, can be written as

$$f(\mathcal{X}(0)) \prod_{t \le Z} \{f(\mathcal{X}(t)|\bar{H}_\mathcal{X}(t))\} \prod_{t < Z} \{f(Y(t)|\mathbf{a}, \bar{H}_\mathcal{X}(t), \mathcal{X}(t))\}^{\delta R(t)}$$

$$\times \prod_{t \le Z} \{1 - E[dN_T(t)|T \ge t, \mathbf{a}, \bar{H}_\mathcal{X}(t), \mathcal{X}(t)]\}^{1 - \delta N_T(t)}$$

$$\times \{E[dN_T(t)|T \ge t, \mathbf{a}, \bar{H}_\mathcal{X}(t), \mathcal{X}(t)]\}^{\delta N_T(t)}$$

$$\times \prod_{t \le Z} \{1 - E[dN_C(t)|\bar{H}_\mathcal{X}(t), \mathcal{X}(t)]\}^{1 - \delta N_C(t)}$$

$$\times \{E[dN_C(t)|\bar{H}_\mathcal{X}(t), \mathcal{X}(t)]\}^{\delta N_C(t)}$$

$$\times \prod_{t < Z} \{1 - E[dR(t)|\bar{H}_\mathcal{X}(t), \mathcal{X}(t)]\}^{1 - \delta R(t)} \{E[dR(t)|\bar{H}_\mathcal{X}(t), \mathcal{X}(t)]\}^{\delta R(t)}.$$

Here $\delta N_T(t) = N_T(t) - N_T(t-)$, $\delta N_C(t) = N_C(t) - N_C(t-)$ and $\delta R(t) = R(t) - R(t-)$.

We further assume that the parameters in $f(\mathcal{X}(0))$, $f(\mathcal{X}(t)|\bar{H}_\mathcal{X}(t))$, $E[dN_C(t)|\bar{H}_\mathcal{X}(t), \mathcal{X}(t)]$ and $E[dR(t)|\bar{H}_\mathcal{X}(t), \mathcal{X}(t)]$ are distinct from the parameters in models (1) and (2). The latter consist of $\boldsymbol{\theta} = (\sigma_y, \text{Vec}(\boldsymbol{\Sigma}_a)^T, \boldsymbol{\beta}^T, \boldsymbol{\gamma}^T, \boldsymbol{\phi}^T)^T$ and $\Lambda(t) = \int_0^t \lambda(s)\,ds$, where $\text{Vec}(\boldsymbol{\Sigma}_a)$ is the vector consisting of all the elements in the upper triangular part of $\boldsymbol{\Sigma}_a$. Then the observed likelihood function of $(\boldsymbol{\theta}, \Lambda)$ from $n$ i.i.d. observations is proportional to

$$\prod_{i=1}^n \int_{\mathbf{a}} \Big\{ (2\pi\sigma_y^2)^{-N_i/2} \exp\{-(\mathbf{Y}_i - \mathbf{X}_i^T\boldsymbol{\beta} - \widetilde{\mathbf{X}}_i^T\mathbf{a})^T(\mathbf{Y}_i - \mathbf{X}_i^T\boldsymbol{\beta} - \widetilde{\mathbf{X}}_i^T\mathbf{a})/2\sigma_y^2\}$$

$$\times \lambda(Z_i)^{\Delta_i} \exp\Big[\Delta_i\{(\boldsymbol{\phi} \circ \widetilde{\mathbf{W}}_i(Z_i))^T\mathbf{a} + \mathbf{W}_i(Z_i)^T\boldsymbol{\gamma}\}$$

(3)
$$- \int_0^{Z_i} e^{(\boldsymbol{\phi} \circ \widetilde{\mathbf{W}}_i(s))^T \mathbf{a} + \mathbf{W}_i(s)^T \boldsymbol{\gamma}} \, d\Lambda(s) \Bigg] \Bigg\}$$

$$\times (2\pi)^{-d_a/2} |\boldsymbol{\Sigma}_a|^{-1/2} \exp\{-\mathbf{a}^T \boldsymbol{\Sigma}_a^{-1} \mathbf{a}/2\} \, d\mathbf{a},$$

where for subject $i$, $\mathbf{Y}_i$ denotes the vector of the observed repeated measurements, $\mathbf{X}_i$ denotes the matrix with each column equal to the observed covariate $\mathbf{X}_i(t)$ at the time of each measurement, $\widetilde{\mathbf{X}}_i$ denotes the matrix with each column equal to the observed covariate $\widetilde{\mathbf{X}}_i(t)$ at the time of each measurement, $N_i$ is the number of the observed repeated measurements and $d_a$ is the dimension of $\mathbf{a}$.

2.2. *Maximum likelihood estimation.* We can obtain the nonparametric maximum likelihood estimates for $\boldsymbol{\theta}$ and $\Lambda$ based on (3). To do that, we let $\Lambda(\cdot)$ be an increasing and right-continuous step-function with jumps only at $Z_i$ for which $\Delta_i = 1$. Moreover, the maximum likelihood estimates maximize a modified object function of (3), which is obtained by replacing $\lambda(Z_i)$ in (3) with $\Lambda\{Z_i\}$, the jump size of $\Lambda(\cdot)$ at $Z_i$. We denote the logarithm of the modified object function as $l_n(\boldsymbol{\theta}, \Lambda)$.

The expectation–maximization (EM) algorithm has been used to calculate the maximum likelihood estimates. In the EM algorithms, the random effects, $\mathbf{a}_i$, $i = 1, \ldots, n$, are treated as missing. Thus, the M-step is to solve the conditional score equations of the complete data given the observed data. Such a conditional expectation is evaluated using the Gaussian–quadrature approximation in the E-step (Chapter 5 of [9]). Since $\Lambda$ can be derived via one-step plug-in in the M-step, such EM algorithms often converge rapidly. However, neither the efficiency nor the convergence of the EM algorithm has been well justified for this context.

We denote the maximum likelihood estimates by $(\hat{\boldsymbol{\theta}}, \hat{\Lambda})$. The profile log-likelihood function for $\boldsymbol{\theta}$ can be used to estimate the asymptotic covariance of $\hat{\boldsymbol{\theta}}$, as given in [18]. In detail, we define the profile log-likelihood function of $\boldsymbol{\theta}$ as $pl_n(\boldsymbol{\theta}) = \max_{\Lambda \in \mathcal{Z}_n} l_n(\boldsymbol{\theta}, \Lambda)$, where $\mathcal{Z}_n$ consists of all the right-continuous step functions only with positive jumps at $Z_i$ for which $\Delta_i = 1$. Then the second-order numerical differences of $pl_n(\boldsymbol{\theta})$ at $\boldsymbol{\theta} = \hat{\boldsymbol{\theta}}$ can be used to approximate the asymptotic variance of $\hat{\boldsymbol{\theta}}$. Especially, for any constant sequence $h_n = O(n^{-1/2})$ and any constant vector $\mathbf{e}$,

$$-(nh_n^2)^{-1}\{pl_n(\hat{\boldsymbol{\theta}} + h_n e) - 2pl_n(\hat{\boldsymbol{\theta}}) + pl_n(\hat{\boldsymbol{\theta}} - h_n e)\}$$

approximates $\mathbf{e}^T \mathbf{I} \mathbf{e}$. Here, $\mathbf{I}$ is the efficient information matrix for $\boldsymbol{\theta}$ and it is also equal to the inverse of the asymptotic covariance of $\sqrt{n}\hat{\boldsymbol{\theta}}$. However, this result is not trivial and will be fully stated and justified in the following sections.

Simulation studies conducted in the past literature [5, 13] have indicated good performance of the maximum likelihood estimates and the proposed variance estimation approach in small samples.

**3. Main results.** In this section we provide the asymptotic properties of $(\hat{\boldsymbol{\theta}}, \hat{\Lambda})$. Theorem 3.1 concerns the consistency of the estimators; Theorem 3.2 gives their asymptotic distribution; the use of the profile log-likelihood function is justified by Theorem 3.3.

In addition to (A.1)–(A.6), some technical assumptions are needed for our main theorems.

(A.7) Recall that $N$ is the number of observed repeated measurements. There exists an integer $n_0$ such that $P(N \leq n_0) = 1$. Moreover, $P(N > d_a | \bar{H}_{\mathcal{X}}(T), \mathcal{X}(T), T) > 0$ with probability one.

(A.8) The maximal right-censoring time is equal to $\tau$.

(A.9) Both $P(\mathbf{X}^T \mathbf{X} \text{ is full rank})$ and $P(\widetilde{\mathbf{X}}^T \widetilde{\mathbf{X}} \text{ is full rank})$ are positive. Additionally, if there exist constant vectors $\mathbf{C}_0$ and $\widetilde{\mathbf{C}}_0$ such that, with positive probability, $\mathbf{W}(t)^T \mathbf{C}_0 = \mu(t), \widetilde{\mathbf{W}}(t) \circ \widetilde{\mathbf{C}}_0 = 0$ for a deterministic function $\mu(t)$ for all $t \in [0, \tau]$, then $\mathbf{C}_0 = 0$, $\widetilde{\mathbf{C}}_0 = 0$ and $\mu(t) \equiv 0$.

(A.10) The true parameter for $\boldsymbol{\theta}$, denoted by $\boldsymbol{\theta}_0 = (\sigma_{0y}, \text{Vec}(\boldsymbol{\Sigma}_{0a})^T, \boldsymbol{\beta}_0^T, \boldsymbol{\gamma}_0^T, \boldsymbol{\phi}_0^T)^T$, satisfies $\|\boldsymbol{\theta}_0\| \leq M_0$, $\sigma_{0y} > M_0^{-1}$, $\min_{\|\mathbf{e}\|=1} \mathbf{e}^T \boldsymbol{\Sigma}_{0a} \mathbf{e} > M_0^{-1}$ for a known positive constant $M_0$.

(A.11) The true hazard rate function, denoted by $\lambda_0(y)$, is bounded and positive in $[0, \tau]$.

REMARK 3.1. Assumption (A.7) stipulates that some subjects have at least $d_a$ repeated measurements. Assumption (A.8) is equivalent to saying that any subject surviving after time $\tau$ is right-censored at $\tau$. Since (A.2), (A.10) and (A.11) imply that, conditional on $\bar{H}_{\mathcal{X}}(\tau)$ and $\mathcal{X}(\tau)$, the probability of a subject surviving after time $\tau$ is at least some positive constant $c_0$, we conclude that $P(C \geq \tau | \bar{H}_{\mathcal{X}}(\tau), \mathcal{X}(\tau)) = P(C = \tau | \bar{H}_{\mathcal{X}}(\tau), \mathcal{X}(\tau)) \geq c_0$ with probability one. The assumptions given in the first half of (A.9) are the same as the identifiability assumption used in a linear mixed effects model, while the second half of (A.9) is used to identify the regression coefficients in the proportional hazards model (1). When $\mathbf{X}(t) = [1, \mathbf{X}]$ and $\widetilde{\mathbf{X}}(t) = \mathbf{W}(t) = \widetilde{\mathbf{W}}(t) = \mathbf{X}$ for time-independent variables $\mathbf{X}$, assumption (A.9) is equivalent to the linear independence of $[1, \mathbf{X}]$ with positive probability. Finally, assumption (A.10) indicates that $\boldsymbol{\theta}_0$ belongs to a known compact set within the domain of $\boldsymbol{\theta}$, denoted by $\Theta$.

We obtain the following theorems.

THEOREM 3.1. *Under assumptions* (A.1)–(A.11), *the maximum likelihood estimator* $(\hat{\boldsymbol{\theta}}, \hat{\Lambda})$ *is strongly consistent under the product metric of the Euclidean norm and the supremum norm on* $[0, \tau]$; *that is,*

$$\|\hat{\boldsymbol{\theta}} - \boldsymbol{\theta}_0\| + \sup_{t \in [0,\tau]} |\hat{\Lambda}(t) - \Lambda_0(t)| \to 0 \qquad a.s.$$

Theorem 3.1 states the strong consistency of the maximum likelihood estimator. Although it is assumed that $\hat{\boldsymbol{\theta}}$ is bounded, we impose no compactness assumption on the estimate $\hat{\Lambda}$. In fact, obtaining the boundedness of $\hat{\Lambda}$ is the key to the proof of Theorem 3.1. The proof will be given in Section 4. Once the result of Theorem 3.1 holds, Theorem 3.2 states the asymptotic properties of the maximum likelihood estimator.

THEOREM 3.2. *Under assumptions* (A.1)–(A.11), $\sqrt{n}(\hat{\boldsymbol{\theta}} - \boldsymbol{\theta}_0, \hat{\Lambda} - \Lambda_0)$ *weakly converges to a Gaussian random element in* $R^d \times l^\infty[0, \tau]$, *where $d$ is the dimension of $\boldsymbol{\theta}$ and $l^\infty[0, \tau]$ is the metric space of all bounded functions in* $[0, \tau]$. *Furthermore,* $\sqrt{n}(\hat{\boldsymbol{\theta}} - \boldsymbol{\theta}_0)$ *weakly converges to a multivariate normal distribution with mean zero and its asymptotic variance attains the semiparametric efficiency bound for* $\boldsymbol{\theta}_0$.

The definition of semiparametric efficiency can be found in Section 3 of [3]. Theorem 3.2 declares that $\hat{\boldsymbol{\theta}}$ is an efficient estimator for $\boldsymbol{\theta}_0$. The proof of Theorem 3.2 is based on the Taylor expansion of the score equations for $\hat{\boldsymbol{\theta}}$ and $\hat{\Lambda}$ around the true parameters $\boldsymbol{\theta}_0$ and $\Lambda_0$. The key to the proof of the theorem is to show that the information operator for $(\boldsymbol{\theta}_0, \Lambda_0)$ is continuously invertible in an appropriate metric space. The proof of Theorem 3.2 will be given in Section 5.

When both Theorems 3.1 and 3.2 are true, we can verify the smooth conditions in Theorem 1 of [18] and show that the profile log-likelihood function, $pl_n(\boldsymbol{\theta})$, approximates a nondegenerate parabolic function around $\hat{\boldsymbol{\theta}}$. Particularly, the inverse of the curvature of the profile log-likelihood function at $\hat{\boldsymbol{\theta}}$ can be used to estimate the asymptotic variance of $\hat{\boldsymbol{\theta}}$. In other words, the following result holds.

THEOREM 3.3. *Under assumptions* (A.1)–(A.11), $2\{pl_n(\hat{\boldsymbol{\theta}}) - pl_n(\boldsymbol{\theta}_0)\}$ *weakly converges to a chi-square distribution with $d$ degrees of freedom and, moreover,*

$$-\frac{pl_n(\hat{\boldsymbol{\theta}} + h_n \mathbf{e}) - 2pl_n(\hat{\boldsymbol{\theta}}) + pl_n(\hat{\boldsymbol{\theta}} - h_n \mathbf{e})}{nh_n^2} \xrightarrow{p} \mathbf{e}^T \mathbf{I} \mathbf{e},$$

*where $h_n = O_p(n^{-1/2})$, $\mathbf{e}$ is any vector in $R^d$ with unit norm, and $\mathbf{I}$ is the efficient information matrix for $\boldsymbol{\theta}_0$.*

REMARK 3.2. Specifically, to estimate the $(s,l)$-element of $\mathbf{I}$, we let $\mathbf{e}_s$ and $\mathbf{e}_l$ be canonical basis elements which have ones at the $s$th coordinate and $l$th coordinate, respectively, and have zeros elsewhere. Then the $(s,l)$-element of $\mathbf{I}$ can be estimated by

$$-(nh_n^2)^{-1}\{pl_n(\hat{\boldsymbol{\theta}} + h_n\mathbf{e}_s + h_n\mathbf{e}_l)$$
$$- pl_n(\hat{\boldsymbol{\theta}} + h_n\mathbf{e}_s - h_n\mathbf{e}_l) - pl_n(\hat{\boldsymbol{\theta}} - h_n\mathbf{e}_s + h_n\mathbf{e}_l) + pl_n(\hat{\boldsymbol{\theta}})\}.$$

**4. Proof of Theorem 3.1.** In this section we prove the consistency result for $(\hat{\boldsymbol{\theta}}, \hat{\Lambda})$. We recall that $l_n(\boldsymbol{\theta}, \Lambda)$ is equal to

$$
\begin{aligned}
&\sum_{i=1}^n \log \int_{\mathbf{a}} \Bigg[ (2\pi\sigma_y^2)^{-N_i/2} \\
&\qquad \times \exp\{-(\mathbf{Y}_i - \mathbf{X}_i^T\boldsymbol{\beta} - \widetilde{\mathbf{X}}_i^T\mathbf{a})^T(\mathbf{Y}_i - \mathbf{X}_i^T\boldsymbol{\beta} - \widetilde{\mathbf{X}}_i^T\mathbf{a})/2\sigma_y^2\} \\
&\qquad \times \Lambda\{Z_i\}^{\Delta_i} \exp\Big\{\Delta_i((\boldsymbol{\phi} \circ \widetilde{\mathbf{W}}_i(Z_i))^T\mathbf{a} + \mathbf{W}_i(Z_i)^T\boldsymbol{\gamma}) \\
&\qquad\qquad - \int_0^{Z_i} e^{(\boldsymbol{\phi} \circ \widetilde{\mathbf{W}}_i(t))^T\mathbf{a} + \mathbf{W}_i(t)^T\boldsymbol{\gamma}}\, d\Lambda(t)\Big\} \\
&\qquad \times (\sqrt{2\pi})^{-d_a}|\boldsymbol{\Sigma}_a|^{-1/2}\exp\{-\mathbf{a}^T\boldsymbol{\Sigma}_a^{-1}\mathbf{a}/2\} \Bigg]\, d\mathbf{a}
\end{aligned}
$$

(4)

and $(\hat{\boldsymbol{\theta}}, \hat{\Lambda})$ maximizes $l_n(\boldsymbol{\theta}, \Lambda)$ over the space $\{(\boldsymbol{\theta}, \Lambda) : \boldsymbol{\theta} \in \Theta, \Lambda \in \mathcal{Z}_n\}$.

The proof of the consistency can be established by verifying the following steps (i)–(iii). One particular remark is that all the following arguments are made and hold for a fixed $\omega$ in the probability space, except for some zero-probability set; thus, all the bounds or constants given below may depend on this $\omega$.

(i) The maximum likelihood estimate $(\hat{\boldsymbol{\theta}}, \hat{\Lambda})$ exists.

(ii) We will show that, with probability one, $\hat{\Lambda}(\tau)$ is bounded as $n$ goes to infinity.

(iii) If (ii) is true, by the Helly selection theorem (cf. page 336 of [4]), we can choose a subsequence of $\hat{\Lambda}$ such that $\hat{\Lambda}$ weakly converges to some right-continuous monotone function $\Lambda^*$ with probability 1; that is, the measure given by $\mu([0,t]) = \hat{\Lambda}(t)$ for $t \in [0,\tau]$ weakly converges to the measure given by $\mu^*([0,t]) = \Lambda^*(t)$. By choosing a sub-subsequence, we can further assume $\hat{\boldsymbol{\theta}} \to \boldsymbol{\theta}^*$. Thus, our third step is to show $\boldsymbol{\theta}^* = \boldsymbol{\theta}_0$ and $\Lambda^* = \Lambda_0$.

Once the above three steps are proved, we conclude that, with probability 1, $\hat{\boldsymbol{\theta}}$ converges to $\boldsymbol{\theta}_0$ and $\hat{\Lambda}$ weakly converges to $\Lambda_0$ in $[0,\tau]$. However, since $\Lambda_0$ is continuous in $[0,\tau]$, the latter can be strengthened to uniform convergence;

that is, $\sup_{t\in[0,\tau]} |\hat{\Lambda}(t) - \Lambda_0(t)| \to 0$ almost surely. Hence, Theorem 3.1 is proved.

The proofs of (i)–(iii) are given in the following.

PROOF OF (i). It suffices to show that the jump size of $\hat{\Lambda}$ at $Z_i$ for which $\Delta_i = 1$ is finite. Since for each $i$ such that $\Delta_i = 1$,

$$\Lambda\{Z_i\}^{\Delta_i} \exp\left\{-\int_0^{Z_i} e^{(\phi \circ \widetilde{\mathbf{W}}_i(t))^T \mathbf{a} + \mathbf{W}_i(t)^T \boldsymbol{\gamma}} \, d\Lambda(t)\right\}$$

$$\leq \exp\{-2((\phi \circ \widetilde{\mathbf{W}}_i(Z_i))^T \mathbf{a} + \mathbf{W}_i(Z_i)^T \boldsymbol{\gamma})\} \Lambda\{Z_i\}^{-1},$$

$l_n(\boldsymbol{\theta}, \Lambda)$ is less than

$$\sum_{i=1}^n \log \int_{\mathbf{a}} \left[ (2\pi\sigma_y^2)^{-N_i/2} \right.$$

$$\times \exp\{-(\mathbf{Y}_i - \mathbf{X}_i^T \boldsymbol{\beta} - \widetilde{\mathbf{X}}_i^T \mathbf{a})^T (\mathbf{Y}_i - \mathbf{X}_i^T \boldsymbol{\beta} - \widetilde{\mathbf{X}}_i^T \mathbf{a})/2\sigma_y^2\}$$

$$\times \exp\{-\Delta_i((\phi \circ \widetilde{\mathbf{W}}_i(Z_i))^T \mathbf{a} + \mathbf{W}_i(Z_i)^T \boldsymbol{\gamma})\}$$

$$\left. \times \Lambda\{Z_i\}^{-\Delta_i} (\sqrt{2\pi})^{-d_a} |\boldsymbol{\Sigma}_a|^{-1/2} \exp\{-\mathbf{a}^T \boldsymbol{\Sigma}_a^{-1} \mathbf{a}/2\} \right] d\mathbf{a}.$$

Thus, if for some $i$ such that $\Delta_i = 1$ and $\Lambda\{Z_i\} \to \infty$, $l_n(\boldsymbol{\theta}, \Lambda) \to -\infty$. We conclude that the jump size of $\Lambda$ must be finite. On the other hand, $\boldsymbol{\theta}$ belongs to a compact set $\Theta$. Then the maximum likelihood estimate $(\hat{\boldsymbol{\theta}}, \hat{\Lambda})$ exists. □

PROOF OF (ii). Define $\hat{\xi} = \log \hat{\Lambda}(\tau)$ and rescale $\hat{\Lambda}$ by the factor $e^{\hat{\xi}}$. We denote the rescaled function as $\widetilde{\Lambda}$; thus $\widetilde{\Lambda}(\tau) = 1$. To prove (ii), it is sufficient to show $\hat{\xi}$ is bounded.

Clearly $\hat{\xi}$ maximizes the log-likelihood function $l_n(\hat{\boldsymbol{\theta}}, \widetilde{\Lambda} e^{\xi})$. After some algebra in expression (4), we obtain that, for any $\Lambda \in \mathcal{Z}_n$, $n^{-1} l_n(\hat{\boldsymbol{\theta}}, \Lambda)$ is equal to

$$-\frac{\sum_{i=1}^n N_i}{n} \log\{\sqrt{2\pi\hat{\sigma}_y^2}\} - \frac{1}{n} \log\{(\sqrt{2\pi})^{d_a} |\hat{\boldsymbol{\Sigma}}_a|^{1/2}\} + \frac{1}{n} \sum_{i=1}^n \log |\mathbf{V}_i|^{1/2}$$

$$+ \frac{1}{n} \sum_{i=1}^n \Delta_i \mathbf{W}_i(Z_i)^T \hat{\boldsymbol{\gamma}} - \frac{1}{n} \sum_{i=1}^n (\mathbf{Y}_i - \mathbf{X}_i^T \hat{\boldsymbol{\beta}})^T (\mathbf{Y}_i - \mathbf{X}_i^T \hat{\boldsymbol{\beta}})/2\hat{\sigma}_y^2$$

$$+ \frac{1}{2n} \sum_{i=1}^n \{(\hat{\boldsymbol{\phi}} \circ \widetilde{\mathbf{W}}_i(Z_i))\Delta_i + \widetilde{\mathbf{X}}_i^T (\mathbf{Y}_i - \mathbf{X}_i^T \hat{\boldsymbol{\beta}})/\hat{\sigma}_y^2\}^T$$

$$\times \mathbf{V}_i^{-1} \{(\hat{\boldsymbol{\phi}} \circ \widetilde{\mathbf{W}}_i(Z_i))\Delta_i + \widetilde{\mathbf{X}}_i^T (\mathbf{Y}_i - \mathbf{X}_i^T \hat{\boldsymbol{\beta}})/\hat{\sigma}_y^2\}$$

$$+ \frac{1}{n}\sum_{i=1}^{n}\left[\Delta_i \log \Lambda\{Z_i\} + \log\left\{\int_{\mathbf{a}} \exp\left\{-\frac{\mathbf{a}^T\mathbf{a}}{2} - \int_0^{Z_i} e^{Q_{1i}(t,\mathbf{a},\hat{\boldsymbol{\theta}})}\,d\Lambda(t)\right\}d\mathbf{a}\right\}\right],$$

where $\mathbf{V}_i = \widetilde{\mathbf{X}}_i^T\widetilde{\mathbf{X}}_i/\hat{\sigma}_y^2 + \hat{\boldsymbol{\Sigma}}_a^{-1}$ and

$$Q_{1i}(t,\mathbf{a},\hat{\boldsymbol{\theta}}) = \{\hat{\boldsymbol{\phi}} \circ \widetilde{\mathbf{W}}_i(t)\}^T \mathbf{V}_i^{-1/2}\mathbf{a} + \mathbf{W}_i(t)^T \hat{\boldsymbol{\gamma}}$$
$$+ \{\hat{\boldsymbol{\phi}} \circ \widetilde{\mathbf{W}}_i(t)\}^T \mathbf{V}_i^{-1}\{(\hat{\boldsymbol{\phi}} \circ \widetilde{\mathbf{W}}_i(Z_i))\Delta_i + \widetilde{\mathbf{X}}_i^T(\mathbf{Y}_i - \mathbf{X}_i^T\hat{\boldsymbol{\beta}})/\hat{\sigma}_y^2\}.$$

Thus, since $0 \leq n^{-1}l_n(\hat{\boldsymbol{\theta}}, e^{\hat{\xi}}\widetilde{\Lambda}) - n^{-1}l_n(\hat{\boldsymbol{\theta}}, \widetilde{\Lambda})$, it follows that

(5)
$$0 \leq \frac{1}{n}\sum_{i=1}^{n}\Delta_i\hat{\xi} + \frac{1}{n}\sum_{i=1}^{n}\log\left\{\int_{\mathbf{a}} \exp\left\{-\frac{\mathbf{a}^T\mathbf{a}}{2} - e^{\hat{\xi}}\int_0^{Z_i} e^{Q_{1i}(t,\mathbf{a},\hat{\boldsymbol{\theta}})}\,d\widetilde{\Lambda}(t)\right\}d\mathbf{a}\right\}$$
$$- \frac{1}{n}\sum_{i=1}^{n}\log\left\{\int_{\mathbf{a}} \exp\left\{-\frac{\mathbf{a}^T\mathbf{a}}{2} - \int_0^{Z_i} e^{Q_{1i}(t,\mathbf{a},\hat{\boldsymbol{\theta}})}\,d\widetilde{\Lambda}(t)\right\}d\mathbf{a}\right\}.$$

According to assumption (A.2), there exist some positive constants $C_1$, $C_2$ and $C_3$ such that $|Q_{1i}(t,\mathbf{a},\hat{\boldsymbol{\theta}})| \leq C_1\|\mathbf{a}\| + C_2\|\mathbf{Y}_i\| + C_3$. If we denote $\mathbf{a}_0$ as the standard multivariate normal distribution, from the concavity of the logarithm function,

$$\log \int_{\mathbf{a}} \exp\left\{-\frac{\mathbf{a}^T\mathbf{a}}{2} - \int_0^{Z_i} e^{Q_{1i}(t,\mathbf{a},\hat{\boldsymbol{\theta}})}\,d\widetilde{\Lambda}(t)\right\}d\mathbf{a}$$
$$= (2\pi)^{d_a/2}\log E_{\mathbf{a}_0}\left[\exp\left\{-\int_0^{Z_i} e^{Q_{1i}(t,\mathbf{a}_0,\hat{\boldsymbol{\theta}})}\,d\widetilde{\Lambda}(t)\right\}\right]$$
$$\geq (2\pi)^{d_a/2}\log E_{\mathbf{a}_0}[\exp\{-e^{C_1\|\mathbf{a}_0\|+C_2\|\mathbf{Y}_i\|+C_3}\}]$$
$$\geq (2\pi)^{d_a/2} E_{\mathbf{a}_0}[-e^{C_1\|\mathbf{a}_0\|+C_2\|\mathbf{Y}_i\|+C_3}]$$
$$= -e^{C_2\|\mathbf{Y}_i\|+C_4},$$

where $C_4$ is another constant. Thus, by the strong law of large numbers and assumption (A.4),

$$-\frac{1}{n}\sum_{i=1}^{n}\log\left\{\int_{\mathbf{a}} \exp\left\{-\frac{\mathbf{a}^T\mathbf{a}}{2} - \int_0^{Z_i} e^{Q_{1i}(t,\mathbf{a},\hat{\boldsymbol{\theta}})}\,d\widetilde{\Lambda}(t)\right\}d\mathbf{a}\right\} \leq \frac{1}{n}\sum_{i=1}^{n}e^{C_2\|\mathbf{Y}_i\|+C_4}$$

can be bounded by some constant $C_5$ from above. Then (5) becomes

$$0 \leq \frac{1}{n}\sum_{i=1}^{n}\Delta_i\hat{\xi}$$
$$+ \frac{1}{n}\sum_{i=1}^{n}\log\left\{\int_{\mathbf{a}} \exp\left\{-\frac{\mathbf{a}^T\mathbf{a}}{2} - e^{\hat{\xi}}\int_0^{Z_i} e^{Q_{1i}(t,\mathbf{a},\hat{\boldsymbol{\theta}})}\,d\widetilde{\Lambda}(t)\right\}d\mathbf{a}\right\} + C_5$$

$$\leq \frac{1}{n}\sum_{i=1}^{n}\Delta_i\xi$$

$$+\frac{1}{n}\sum_{i=1}^{n}I(Z_i=\tau)\log\left\{\int_{\mathbf{a}}\exp\left\{-\frac{\mathbf{a}^T\mathbf{a}}{2}-e^{\hat{\xi}}\int_0^{\tau}e^{Q_{1i}(t,\mathbf{a},\hat{\boldsymbol{\theta}})}\,d\widetilde{\Lambda}(t)\right\}d\mathbf{a}\right\}$$

(6)

$$+\frac{1}{n}\sum_{i=1}^{n}I(Z_i\neq\tau)\log\left\{\int_{\mathbf{a}}\exp\left\{-\frac{\mathbf{a}^T\mathbf{a}}{2}\right\}d\mathbf{a}\right\}+C_5$$

$$\leq \frac{1}{n}\sum_{i=1}^{n}\Delta_i\xi$$

$$+\frac{1}{n}\sum_{i=1}^{n}I(Z_i=\tau)\log\left\{\int_{\mathbf{a}}\exp\left\{-\frac{\mathbf{a}^T\mathbf{a}}{2}-e^{\hat{\xi}}\int_0^{\tau}e^{Q_{1i}(t,\mathbf{a},\hat{\boldsymbol{\theta}})}\,d\widetilde{\Lambda}(t)\right\}d\mathbf{a}\right\}$$

$$+C_6,$$

where $C_6$ is a constant.

On the other hand, since for any $\Gamma\geq 0$ and $x>0$, $\Gamma\log(1+\frac{x}{\Gamma})\leq \Gamma\cdot\frac{x}{\Gamma}=x$, we have that $e^{-x}\leq (1+\frac{x}{\Gamma})^{-\Gamma}$. Therefore,

$$\exp\left\{-\frac{\mathbf{a}^T\mathbf{a}}{2}-e^{\hat{\xi}}\int_0^{\tau}e^{Q_{1i}(t,\mathbf{a},\hat{\boldsymbol{\theta}})}\,d\widetilde{\Lambda}(t)\right\}$$

$$\leq \exp\left\{-\frac{\mathbf{a}^T\mathbf{a}}{2}\right\}\left\{1+e^{\hat{\xi}}\int_0^{\tau}e^{Q_{1i}(t,\mathbf{a},\hat{\boldsymbol{\theta}})}\,d\widetilde{\Lambda}(t)/\Gamma\right\}^{-\Gamma}$$

$$\leq \int_{\mathbf{a}}\Gamma^{\Gamma}\exp\left\{-\frac{\mathbf{a}^T\mathbf{a}}{2}\right\}\left\{e^{\hat{\xi}}\int_0^{\tau}e^{Q_{1i}(t,\mathbf{a},\hat{\boldsymbol{\theta}})}\,d\widetilde{\Lambda}(t)\right\}^{-\Gamma}d\mathbf{a}$$

$$\leq \int_{\mathbf{a}}\Gamma^{\Gamma}\exp\left\{-\frac{\mathbf{a}^T\mathbf{a}}{2}-\Gamma\hat{\xi}\right\}\left\{\int_0^{\tau}e^{Q_{1i}(t,\mathbf{a},\hat{\boldsymbol{\theta}})}\,d\widetilde{\Lambda}(t)\right\}^{-\Gamma}d\mathbf{a}.$$

Since $Q_{1i}(t,\mathbf{a},\hat{\boldsymbol{\theta}})\geq -C_1\|\mathbf{a}\|-C_2\|\mathbf{Y}_i\|-C_3$, (6) gives that

$$0\leq C_6+\frac{1}{n}\sum_{i=1}^{n}\Delta_i\hat{\xi}$$

$$+\frac{1}{n}\sum_{i=1}^{n}I(Z_i=\tau)\log\left\{e^{-\Gamma\hat{\xi}}\Gamma^{\Gamma}\right.$$

$$\left.\times\int_{\mathbf{a}}\exp\left\{-\frac{\mathbf{a}^T\mathbf{a}}{2}+C_1\Gamma\|\mathbf{a}\|+C_2\Gamma\|\mathbf{Y}_i\|+C_3\Gamma\right\}d\mathbf{a}\right\}$$

$$\leq C_6+\frac{1}{n}\sum_{i=1}^{n}\Delta_i\hat{\xi}-\frac{\Gamma}{n}\sum_{i=1}^{n}I(Z_i=\tau)\hat{\xi}+C_7(\Gamma),$$

where $C_7(\Gamma)$ is a deterministic function of $\Gamma$. By the strong law of large numbers, $n^{-1}\sum_{i=1}^n I(Z_i = \tau) \to P(Z = \tau) > 0$. Then we can choose $\Gamma$ large enough such that $n^{-1}\sum_{i=1}^n \Delta_i \leq (2n)^{-1}\Gamma \sum_{i=1}^n I(Z_i = \tau)$. Thus, we obtain that

$$0 \leq C_7(\Gamma) + C_6 - \frac{\Gamma}{2n}\sum_{i=1}^n I(Z_i = \tau)\hat{\xi}.$$

In other words, if we let $B_0 = \exp\{4(C_6 + C_7(\Gamma))/\Gamma P(Z = \tau)\}$, we conclude that $\hat{\Lambda}(\tau) \leq B_0$. Note that the above arguments hold for every sample in the probability space except a set with zero probability. Thus, we have shown that, with probability 1, $\hat{\Lambda}(\tau)$ is bounded for any sample size $n$. □

PROOF OF (iii). In this step we will show that if $\hat{\boldsymbol{\theta}} \to \boldsymbol{\theta}^*$ and $\hat{\Lambda}$ weakly converges to $\Lambda^*$ with probability one, then $\boldsymbol{\theta}^* = \boldsymbol{\theta}_0$ and $\Lambda^* = \Lambda_0$. For convenience, we use $\mathbf{O}$ to abbreviate the observed statistics $(\mathbf{Y}, \mathbf{X}, \widetilde{\mathbf{X}}, Z, \Delta, N)$ and $\{\mathbf{W}(s), \widetilde{\mathbf{W}}(s), 0 \leq s \leq Z\}$ and denote $G(\mathbf{a}, \mathbf{O}; \boldsymbol{\theta}, \Lambda)$ as

$$(2\pi\sigma_y^2)^{-N/2}\{(2\pi)^{d_a}|\boldsymbol{\Sigma}_a|\}^{-1/2}$$
$$\times \exp\bigg\{-\frac{(\mathbf{Y} - \mathbf{X}^T\boldsymbol{\beta} - \widetilde{\mathbf{X}}^T\mathbf{a})^T(\mathbf{Y} - \mathbf{X}^T\boldsymbol{\beta} - \widetilde{\mathbf{X}}^T\mathbf{a})}{2\sigma_y^2} - \frac{\mathbf{a}^T\boldsymbol{\Sigma}_a^{-1}\mathbf{a}}{2}$$
$$+ \Delta((\boldsymbol{\phi}\circ\widetilde{\mathbf{W}}(Z))^T\mathbf{a} + \mathbf{W}(Z)^T\boldsymbol{\gamma}) - \int_0^Z e^{(\boldsymbol{\phi}\circ\widetilde{\mathbf{W}}(t))^T\mathbf{a} + \mathbf{W}(t)^T\boldsymbol{\gamma}}d\Lambda(t)\bigg\}.$$

Moreover, we define

$$Q(z, \mathbf{O}; \boldsymbol{\theta}, \Lambda) = \frac{\int_{\mathbf{a}} G(\mathbf{a}, \mathbf{O}; \boldsymbol{\theta}, \Lambda)\exp\{(\boldsymbol{\phi}\circ\widetilde{\mathbf{W}}(z))^T\mathbf{a} + \mathbf{W}(z)^T\boldsymbol{\gamma}\}d\mathbf{a}}{\int_{\mathbf{a}} G(\mathbf{a}, \mathbf{O}; \boldsymbol{\theta}, \Lambda)d\mathbf{a}},$$

and for any measurable function $f(\mathbf{O})$, we use operator notation to define $\mathbf{P}_n f = n^{-1}\sum_{i=1}^n f(\mathbf{O}_i)$ and $\mathbf{P}f = \int f\,d\mathbf{P} = E[f(\mathbf{O})]$. Thus, $\mathbf{P}_n$ is the empirical measure from $n$ i.i.d. observations and $\sqrt{n}(\mathbf{P}_n - \mathbf{P})$ is the empirical process based on these observations (cf. Section 2.1 of [24]). We also define a class $\mathcal{F} = \{Q(z, \mathbf{O}; \boldsymbol{\theta}, \Lambda) : z \in [0, \tau], \boldsymbol{\theta} \in \Theta, \Lambda \in \mathcal{Z}, \Lambda(0) = 0, \Lambda(\tau) \leq B_0\}$, where $B_0$ is the constant given in step (ii) and $\mathcal{Z}$ contains all nondecreasing functions in $[0, \tau]$. According to Appendix A.1, $\mathcal{F}$ is P-Donsker.

We start to prove (iii). Since $(\hat{\boldsymbol{\theta}}, \hat{\Lambda})$ maximizes the function $l_n(\boldsymbol{\theta}, \Lambda)$, where $\Lambda$ is any step function with jumps only at $Z_i$ for which $\Delta_i = 1$, after differentiating $l_n(\boldsymbol{\theta}, \Lambda)$ with respect to $\Lambda\{Z_i\}$, we obtain that $\hat{\Lambda}$ satisfies the equation

$$\hat{\Lambda}\{Z_k\} = \Delta_k/[n\mathbf{P}_n\{I(Z \geq z)Q(z, \mathbf{O}; \hat{\boldsymbol{\theta}}, \hat{\Lambda})\}]|_{z=Z_k}.$$

Thus, imitating the above equation, we can construct another function, denoted by $\bar{\Lambda}$, such that $\bar{\Lambda}$ is also a step function with jumps only at observed $Z_k$ and the jump size is given by

$$\bar{\Lambda}\{Z_k\} = \Delta_k/[n\mathbf{P}_n\{I(Z \geq z)Q(z, \mathbf{O}; \boldsymbol{\theta}_0, \Lambda_0)\}]|_{z=Z_k}.$$

Equivalently,

$$\bar{\Lambda}(t) = \frac{1}{n}\sum_{k=1}^{n}\frac{I(Z_k \leq t)\Delta_k}{\mathbf{P}_n\{I(Z \geq z)Q(z, \mathbf{O}; \boldsymbol{\theta}_0, \Lambda_0)\}|_{z=Z_k}}.$$

We claim $\bar{\Lambda}(t)$ uniformly converges to $\Lambda_0(t)$ in $[0,\tau]$. To prove the claim, we note that

$$\sup_{t\in[0,\tau]}\left|\bar{\Lambda}(t) - E\left[\frac{I(Z \leq t)\Delta}{\mathbf{P}\{I(Z \geq z)Q(z, \mathbf{O}; \boldsymbol{\theta}_0, \Lambda_0)\}|_{z=Z}}\right]\right|$$

$$\leq \sup_{t\in[0,\tau]}\left|\frac{1}{n}\sum_{k=1}^{n}I(Z_k \leq t)\Delta_k\left[\frac{1}{\mathbf{P}_n\{I(Z \geq z)Q(z, \mathbf{O}; \boldsymbol{\theta}_0, \Lambda_0)\}}\right.\right.$$

$$\left.\left. - \frac{1}{\mathbf{P}\{I(Z \geq z)Q(z, \mathbf{O}; \boldsymbol{\theta}_0, \Lambda_0)\}}\right]\right|_{z=Z_k}$$

$$+ \sup_{t\in[0,\tau]}\left|(\mathbf{P}_n - \mathbf{P})\left[\frac{I(Z \leq t)\Delta}{\mathbf{P}\{I(Z \geq z)Q(z, \mathbf{O}; \boldsymbol{\theta}_0, \Lambda_0)\}|_{z=Z}}\right]\right|$$

$$\leq \sup_{z\in[0,\tau]}\left|\frac{1}{\mathbf{P}_n\{I(Z \geq z)Q(z, \mathbf{O}; \boldsymbol{\theta}_0, \Lambda_0)\}} - \frac{1}{\mathbf{P}\{I(Z \geq z)Q(z, \mathbf{O}; \boldsymbol{\theta}_0, \Lambda_0)\}}\right|$$

$$+ \sup_{t\in[0,\tau]}\left|(\mathbf{P}_n - \mathbf{P})\left[\frac{I(Z \leq t)\Delta}{\mathbf{P}\{I(Z \geq z)Q(z, \mathbf{O}; \boldsymbol{\theta}_0, \Lambda_0)\}|_{z=Z}}\right]\right|.$$

According to Appendix A.1, $\{Q(z, \mathbf{O}; \boldsymbol{\theta}_0, \Lambda_0): z \in [0,\tau]\}$ is a bounded and Glivenko–Cantelli class. Since $\{I(Z \geq z): z \in [0,\tau]\}$ is also a Glivenko–Cantelli class and the functional $(f,g) \mapsto fg$ for any bounded two functions $f$ and $g$ is Lipschitz continuous, $\{I(Z \geq z)Q(z, \mathbf{O}; \boldsymbol{\theta}_0, \Lambda_0): z \in [0,\tau]\}$ is a Glivenko–Cantelli class. Then we obtain that $\sup_{z\in[0,\tau]}|\mathbf{P}_n\{I(Z \geq z)Q(z, \mathbf{O}; \boldsymbol{\theta}_0, \Lambda_0)\} - \mathbf{P}\{I(Z \geq z)Q(z, \mathbf{O}; \boldsymbol{\theta}_0, \Lambda_0)\}|$ converges to 0. Moreover, from Appendix A.1, $\mathbf{P}\{I(Z \geq z)Q(z, \mathbf{O}; \boldsymbol{\theta}_0, \Lambda_0)\}$ is larger than $\mathbf{P}\{I(Z \geq \tau)\exp\{-C_8 - C_9\|\mathbf{Y}\|\}\}$ for two constants $C_8$ and $C_9$, so is bounded from below. Thus, the first term on the right-hand side of the above inequality tends to zero. Additionally, since the class $\{I(Z \leq t)/\mathbf{P}\{I(Z \geq z)Q(z, \mathbf{O}; \boldsymbol{\theta}_0, \Lambda_0)\}|_{z=Z}: t \in [0,\tau]\}$ is also a Glivenko–Cantelli class, the second term on the right-hand side of the above inequality vanishes as $n$ goes to infinity. Therefore, we conclude that $\bar{\Lambda}(t)$ uniformly converges to

$$E\left[\frac{I(Z \leq t)\Delta}{\mathbf{P}\{I(Z \geq z)Q(z, \mathbf{O}; \boldsymbol{\theta}_0, \Lambda_0)\}|_{z=Z}}\right].$$

It is easy to verify that this limit is equal to $\Lambda_0(t)$. Thus, our claim that $\bar{\Lambda}$ uniformly converges to $\Lambda_0$ in $[0,\tau]$ holds.

From the construction of $\bar{\Lambda}$, we obtain that

$$\hat{\Lambda}(t) = \int_0^t \frac{\mathbf{P}_n\{I(Z \geq z)Q(z,\mathbf{O};\boldsymbol{\theta}_0,\Lambda_0)\}}{\mathbf{P}_n\{I(Z \geq z)Q(z,\mathbf{O};\hat{\boldsymbol{\theta}},\hat{\Lambda})\}} \, d\bar{\Lambda}(z). \tag{7}$$

$\hat{\Lambda}(t)$ is absolutely continuous with respect to $\bar{\Lambda}(t)$. On the other hand, since $\{I(Z \geq z) : z \in [0,\tau]\}$ and $\mathcal{F}$ are both Glivenko–Cantelli classes, $\{I(Z \geq z)Q(z,\mathbf{O};\boldsymbol{\theta},\Lambda) : z \in [0,\tau], \boldsymbol{\theta} \in \Theta, \Lambda \in \mathcal{Z}, \Lambda(\tau) \leq B_0\}$ is also a Glivenko–Cantelli class. Thus,

$$\sup_{z \in [0,\tau]} |(\mathbf{P}_n - \mathbf{P})\{I(Z \geq z)Q(z,\mathbf{O};\hat{\boldsymbol{\theta}},\hat{\Lambda})\}|$$
$$+ \sup_{z \in [0,\tau]} |(\mathbf{P}_n - \mathbf{P})\{I(Z \geq z)Q(z,\mathbf{O};\boldsymbol{\theta}_0,\Lambda_0)\}| \to 0 \quad \text{a.s.}$$

On the other hand, using the bounded convergence theorem and the fact that $\hat{\boldsymbol{\theta}}$ converges to $\boldsymbol{\theta}^*$ and $\hat{\Lambda}$ weakly converges to $\Lambda^*$, $\mathbf{P}\{I(Z \geq z)Q(z,\mathbf{O};\hat{\boldsymbol{\theta}},\hat{\Lambda})\}$ converges to $\mathbf{P}\{I(Z \geq z)Q(z,\mathbf{O};\boldsymbol{\theta}^*,\Lambda^*)\}$ for each $z$; moreover, it is straightforward to check that the derivative of $\mathbf{P}\{I(Z \geq z)Q(z,\mathbf{O};\hat{\boldsymbol{\theta}},\hat{\Lambda})\}$ with respect to $z$ is uniformly bounded, so $\mathbf{P}\{I(Z \geq z)Q(z,\mathbf{O};\hat{\boldsymbol{\theta}},\hat{\Lambda})\}$ is equi-continuous with respect to $z$. Thus, by the Arzela–Ascoli theorem (page 245 of [21]), uniformly in $z \in [0,\tau]$,

$$\mathbf{P}\{I(Z \geq z)Q(z,\mathbf{O};\hat{\boldsymbol{\theta}},\hat{\Lambda})\} \to \mathbf{P}\{I(Z \geq z)Q(z,\mathbf{O};\boldsymbol{\theta}^*,\Lambda^*)\}.$$

Then it holds that, uniformly in $z \in [0,\tau]$,

$$\frac{\hat{\Lambda}\{z\}}{\bar{\Lambda}\{z\}} = \frac{\mathbf{P}_n\{I(Z \geq z)Q(z,\mathbf{O};\boldsymbol{\theta}_0,\Lambda_0)\}}{\mathbf{P}_n\{I(Z \geq z)Q(z,\mathbf{O};\hat{\boldsymbol{\theta}},\hat{\Lambda})\}} \to \frac{\mathbf{P}\{I(Z \geq z)Q(z,\mathbf{O};\boldsymbol{\theta}_0,\Lambda_0)\}}{\mathbf{P}\{I(Z \geq z)Q(z,\mathbf{O};\boldsymbol{\theta}^*,\Lambda^*)\}}. \tag{8}$$

After taking limits on both sides of (7), we obtain that

$$\Lambda^*(t) = \int_0^t \frac{\mathbf{P}\{I(Z \geq z)Q(z,\mathbf{O};\boldsymbol{\theta}_0,\Lambda_0)\}}{\mathbf{P}\{I(Z \geq z)Q(z,\mathbf{O};\boldsymbol{\theta}^*,\Lambda^*)\}} \, d\Lambda_0(z).$$

Therefore, since $\Lambda_0(t)$ is differentiable with respect to the Lebesgue measure, so is $\Lambda^*(t)$ and we denote $\lambda^*(t)$ as the derivative of $\Lambda^*(t)$. Additionally, from (8) we note that $\hat{\Lambda}\{Z\}/\bar{\Lambda}\{Z\}$ uniformly converges to $d\Lambda^*(Z)/d\Lambda_0(Z) = \lambda^*(Z)/\lambda_0(Z)$. A second conclusion is that $\hat{\Lambda}$ uniformly converges to $\Lambda^*$ since $\Lambda^*$ is continuous.

On the other hand,

$$n^{-1}l_n(\hat{\boldsymbol{\theta}},\hat{\Lambda}) - n^{-1}l_n(\boldsymbol{\theta}_0,\bar{\Lambda})$$
$$= \mathbf{P}_n\left[\Delta \log \frac{\hat{\Lambda}\{Z\}}{\bar{\Lambda}\{Z\}}\right] + \mathbf{P}_n\left[\log \frac{\int_{\mathbf{a}} G(\mathbf{a},\mathbf{O};\hat{\boldsymbol{\theta}},\hat{\Lambda}) \, d\mathbf{a}}{\int_{\mathbf{a}} G(\mathbf{a},\mathbf{O};\boldsymbol{\theta}_0,\bar{\Lambda}) \, d\mathbf{a}}\right]$$
$$\geq 0.$$

Using the result of Appendix A.1 and similar arguments as above, we can verify that $\log[\int_{\mathbf{a}} G(\mathbf{a}, \mathbf{O}; \hat{\boldsymbol{\theta}}, \hat{\Lambda}) \, d\mathbf{a} / \int_{\mathbf{a}} G(\mathbf{a}, \mathbf{O}; \boldsymbol{\theta}_0, \bar{\Lambda}) \, d\mathbf{a}]$ belongs to a Glivenko–Cantelli class and

$$\mathbf{P}\bigg[\log \frac{\int_{\mathbf{a}} G(\mathbf{a}, \mathbf{O}; \hat{\boldsymbol{\theta}}, \hat{\Lambda}) \, d\mathbf{a}}{\int_{\mathbf{a}} G(\mathbf{a}, \mathbf{O}; \boldsymbol{\theta}_0, \bar{\Lambda}) \, d\mathbf{a}}\bigg] \to \mathbf{P}\bigg[\log \frac{\int_{\mathbf{a}} G(\mathbf{a}, \mathbf{O}; \boldsymbol{\theta}^*, \Lambda^*) \, d\mathbf{a}}{\int_{\mathbf{a}} G(\mathbf{a}, \mathbf{O}; \boldsymbol{\theta}_0, \Lambda_0) \, d\mathbf{a}}\bigg].$$

Since $\hat{\Lambda}\{Z\}/\bar{\Lambda}\{Z\}$ uniformly converges to $\lambda^*(Z)/\lambda_0(Z)$, we obtain that

$$\mathbf{P}\bigg[\log\bigg\{\frac{\lambda^*(Z)^\Delta \int_{\mathbf{a}} G(\mathbf{a}, \mathbf{O}; \boldsymbol{\theta}^*, \Lambda^*) \, d\mathbf{a}}{\lambda_0(Z)^\Delta \int_{\mathbf{a}} G(\mathbf{a}, \mathbf{O}; \boldsymbol{\theta}_0, \Lambda_0) \, d\mathbf{a}}\bigg\}\bigg] \geq 0.$$

However, the left-hand side of the inequality is the negative Kullback–Leibler information. Then it immediately follows that, with probability one,

$$(9) \qquad \lambda^*(Z)^\Delta \int_{\mathbf{a}} G(\mathbf{a}, \mathbf{O}; \boldsymbol{\theta}^*, \Lambda^*) \, d\mathbf{a} = \lambda_0(Z)^\Delta \int_{\mathbf{a}} G(\mathbf{a}, \mathbf{O}; \boldsymbol{\theta}_0, \Lambda_0) \, d\mathbf{a}.$$

Our proof is completed if we can show $\boldsymbol{\theta}^* = \boldsymbol{\theta}_0$ and $\Lambda^* = \Lambda_0$ from (9). Since (9) holds with probability one, (9) holds for any $(Z, \Delta = 1)$ and $(Z = \tau, \Delta = 0)$, but may not hold for $(Z, \Delta = 0)$ when $Z \in (0, \tau)$. However, we can show that (9) is also true for $(Z, \Delta = 0)$ when $Z \in (0, \tau)$. To see that, treating both sides of (9) as functions of $Z$ we integrate these functions over the interval $(Z, \tau)$ to obtain

$$\int_{\mathbf{a}} G(\mathbf{a}, \mathbf{O}; \boldsymbol{\theta}^*, \Lambda^*) \, d\mathbf{a}|_{\Delta=0, Z=\tau} - \int_{\mathbf{a}} G(\mathbf{a}, \mathbf{O}; \boldsymbol{\theta}^*, \Lambda^*) \, d\mathbf{a}|_{\Delta=0, Z=Z}$$
$$= \int_{\mathbf{a}} G(\mathbf{a}, \mathbf{O}; \boldsymbol{\theta}_0, \Lambda_0) \, d\mathbf{a}|_{\Delta=0, Z=\tau} - \int_{\mathbf{a}} G(\mathbf{a}, \mathbf{O}; \boldsymbol{\theta}_0, \Lambda_0) \, d\mathbf{a}|_{\Delta=0, Z=Z}.$$

After comparing this equality with another equality, which is given by (9) at $\Delta = 0$ and $Z = \tau$, we obtain

$$\int_{\mathbf{a}} G(\mathbf{a}, \mathbf{O}; \boldsymbol{\theta}^*, \Lambda^*) \, d\mathbf{a}|_{\Delta=0} = \int_{\mathbf{a}} G(\mathbf{a}, \mathbf{O}; \boldsymbol{\theta}_0, \Lambda_0) \, d\mathbf{a}|_{\Delta=0};$$

that is, (9) also holds for any $Z$ and $\Delta = 0$.

Thus, we let $\Delta = 0$ and $Z = 0$ in (9). After integrating over $\mathbf{a}$, we have that with probability one,

$$\frac{1}{\sigma_y^{*N}} \frac{|\boldsymbol{\Sigma}_a^{*-1} + \widetilde{\mathbf{X}}^T\widetilde{\mathbf{X}}/\sigma_y^{*2}|^{-1/2}}{|\boldsymbol{\Sigma}_a^*|^{1/2}}$$
$$\times \exp\bigg\{-\frac{(\mathbf{Y} - \mathbf{X}^T\boldsymbol{\beta}^*)^T(\mathbf{Y} - \mathbf{X}^T\boldsymbol{\beta}^*)}{2\sigma_y^{*2}}$$
$$+ \frac{1}{2\sigma_y^{*4}}(\mathbf{Y} - \mathbf{X}^T\boldsymbol{\beta}^*)^T\widetilde{\mathbf{X}}\bigg(\boldsymbol{\Sigma}_a^{*-1} + \frac{\widetilde{\mathbf{X}}^T\widetilde{\mathbf{X}}}{\sigma_y^{*2}}\bigg)^{-1}\widetilde{\mathbf{X}}^T(\mathbf{Y} - \mathbf{X}^T\boldsymbol{\beta}^*)\bigg\}$$

(10)
$$= \frac{1}{\sigma_{0y}^N} \frac{|\boldsymbol{\Sigma}_{0a}^{-1} + \widetilde{\mathbf{X}}^T\widetilde{\mathbf{X}}/\sigma_{0y}^2|^{-1/2}}{|\boldsymbol{\Sigma}_{0a}|^{1/2}}$$
$$\times \exp\left\{ -\frac{(\mathbf{Y}-\mathbf{X}^T\boldsymbol{\beta}_0)^T(\mathbf{Y}-\mathbf{X}^T\boldsymbol{\beta}_0)}{2\sigma_{0y}^2} \right.$$
$$\left. + \frac{1}{2\sigma_{0y}^4}(\mathbf{Y}-\mathbf{X}^T\boldsymbol{\beta}_0)^T\widetilde{\mathbf{X}}\left(\boldsymbol{\Sigma}_{0a}^{-1} + \frac{\widetilde{\mathbf{X}}^T\widetilde{\mathbf{X}}}{\sigma_{0y}^2}\right)^{-1}\widetilde{\mathbf{X}}^T(\mathbf{Y}-\mathbf{X}^T\boldsymbol{\beta}_0) \right\}.$$

By comparing the coefficients of $\mathbf{Y}\mathbf{Y}^T$, $\mathbf{Y}$ and the constant term in the exponential parts, we obtain

(11)
$$\frac{1}{\sigma_y^{*4}} \widetilde{\mathbf{X}} \left( \boldsymbol{\Sigma}_a^{*-1} + \frac{\widetilde{\mathbf{X}}^T\widetilde{\mathbf{X}}}{\sigma_y^{*2}} \right)^{-1} \widetilde{\mathbf{X}}^T - \frac{1}{\sigma_y^{*2}}$$
$$= \frac{1}{\sigma_{0y}^4} \widetilde{\mathbf{X}} \left( \boldsymbol{\Sigma}_{0a}^{-1} + \frac{\widetilde{\mathbf{X}}^T\widetilde{\mathbf{X}}}{\sigma_{0y}^2} \right)^{-1} \widetilde{\mathbf{X}}^T - \frac{1}{\sigma_{0y}^2},$$

(12)
$$\mathbf{X}^T\boldsymbol{\beta}^* = \mathbf{X}^T\boldsymbol{\beta}_0$$

and

(13)
$$\frac{1}{\sigma_y^{*N}} \frac{|\boldsymbol{\Sigma}_a^{*-1} + \widetilde{\mathbf{X}}^T\widetilde{\mathbf{X}}/\sigma_y^{*2}|^{-1/2}}{|\boldsymbol{\Sigma}_a^*|^{1/2}} = \frac{1}{\sigma_{0y}^N} \frac{|\boldsymbol{\Sigma}_{0a}^{-1} + \widetilde{\mathbf{X}}^T\widetilde{\mathbf{X}}/\sigma_{0y}^2|^{-1/2}}{|\boldsymbol{\Sigma}_{0a}|^{1/2}}.$$

By (A.9), (12) gives $\boldsymbol{\beta}^* = \boldsymbol{\beta}_0$. We multiply both sides of (11) by $\widetilde{\mathbf{X}}^T$ from the left and by $\widetilde{\mathbf{X}}$ from the right. According to assumption (A.9), it holds that, with positive probability,

$$\frac{1}{\sigma_y^{*4}} \left( \boldsymbol{\Sigma}_a^{*-1} + \frac{\widetilde{\mathbf{X}}^T\widetilde{\mathbf{X}}}{\sigma_y^{*2}} \right)^{-1} \widetilde{\mathbf{X}}^T\widetilde{\mathbf{X}} - \frac{1}{\sigma_y^{*2}} = \frac{1}{\sigma_{0y}^4} \left( \boldsymbol{\Sigma}_{0a}^{-1} + \frac{\widetilde{\mathbf{X}}^T\widetilde{\mathbf{X}}}{\sigma_{0y}^2} \right)^{-1} \widetilde{\mathbf{X}}^T\widetilde{\mathbf{X}} - \frac{1}{\sigma_{0y}^2}.$$

Thus,
$$\sigma_y^{*2} + \boldsymbol{\Sigma}_a^*\widetilde{\mathbf{X}}^T\widetilde{\mathbf{X}} = \sigma_{0y}^2 + \boldsymbol{\Sigma}_{0a}\widetilde{\mathbf{X}}^T\widetilde{\mathbf{X}}.$$

Combining this result with (13) and by assumption (A.7), where $P(N > d_a|\bar{H}_\mathcal{X}(\tau), \mathcal{X}(\tau)) > 0$, we have that $\sigma_y^* = \sigma_{0y}$. This further gives $\boldsymbol{\Sigma}_a^* = \boldsymbol{\Sigma}_{0a}$.

Next, to show that $\boldsymbol{\phi}^* = \boldsymbol{\phi}_0, \boldsymbol{\gamma}^* = \boldsymbol{\gamma}_0$ and $\Lambda^* = \Lambda_0$, we let $\Delta = 0$ in (9) and notice that (9) can be written as

$$E_\mathbf{a}\left[\exp\left\{-\int_0^Z e^{(\boldsymbol{\phi}^* \circ \widetilde{\mathbf{W}}(t))^T \mathbf{a} + \mathbf{W}(t)^T \boldsymbol{\gamma}^*} \, d\Lambda^*(t)\right\}\right]$$
$$= E_\mathbf{a}\left[\exp\left\{-\int_0^Z e^{(\boldsymbol{\phi}_0 \circ \widetilde{\mathbf{W}}(t))^T \mathbf{a} + \mathbf{W}(t)^T \boldsymbol{\gamma}_0} \, d\Lambda_0(t)\right\}\right],$$

where **a** follows a normal distribution with mean $\{\widetilde{\mathbf{X}}^T\widetilde{\mathbf{X}}/\sigma_{0y}^2 + \boldsymbol{\Sigma}_{0a}^{-1}\}\{(\mathbf{Y} - \mathbf{X}^T\boldsymbol{\beta}_0)^T\widetilde{\mathbf{X}}/\sigma_{0y}^2\}^T$ and covariance $[\widetilde{\mathbf{X}}^T\widetilde{\mathbf{X}}/\sigma_{0y}^2 + \boldsymbol{\Sigma}_{0a}^{-1}]^{-1}$. However, for any fixed $\widetilde{\mathbf{X}}$ and $\mathbf{X}$, treating $\widetilde{\mathbf{X}}^T\mathbf{Y}$ as parameters in this normal family, **a** is the complete statistic for $\widetilde{\mathbf{X}}^T\mathbf{Y}$. Therefore,

$$\exp\left\{-\int_0^Z e^{(\boldsymbol{\phi}^*\circ\widetilde{\mathbf{W}}(t))^T\mathbf{a}+\mathbf{W}(t)^T\boldsymbol{\gamma}^*}\,d\Lambda^*(t)\right\}$$
$$=\exp\left\{-\int_0^Z e^{(\boldsymbol{\phi}_0\circ\widetilde{\mathbf{W}}(t))^T\mathbf{a}+\mathbf{W}(t)^T\boldsymbol{\gamma}_0}\,d\Lambda_0(t)\right\}.$$

Equivalently,

$$e^{(\boldsymbol{\phi}^*\circ\widetilde{\mathbf{W}}(t))^T\mathbf{a}+\mathbf{W}(t)^T\boldsymbol{\gamma}^*}\lambda^*(t) = e^{(\boldsymbol{\phi}_0\circ\widetilde{\mathbf{W}}(t))^T\mathbf{a}+\mathbf{W}(t)^T\boldsymbol{\gamma}_0}\lambda_0(t).$$

According to assumptions (A.9) and (A.11), $\boldsymbol{\phi}^* = \boldsymbol{\phi}_0, \boldsymbol{\gamma}^* = \boldsymbol{\gamma}_0$ and $\Lambda^* = \Lambda_0$. □

**5. Proof of Theorem 3.2.** The asymptotic properties for the estimators $(\hat{\boldsymbol{\theta}}, \hat{\Lambda})$ follow if we can verify the conditions in [24], Theorem 3.3.1. For completeness, we state this theorem below (the version from Appendix A of [19]).

THEOREM 5.1 (Theorem 3.3.1 in [24]). *Let $S_n$ and $S$ be random maps and a fixed map, respectively, from $\psi$ to a Banach space such that:*

(a) $\sqrt{n}(S_n - S)(\hat{\psi}_n) - \sqrt{n}(S_n - S)(\psi_0) = o_P^*(1 + \sqrt{n}\|\hat{\psi}_n - \psi_0\|).$
(b) *The sequence $\sqrt{n}(S_n - S)(\psi_0)$ converges in distribution to a tight random element $\mathbf{Z}$.*
(c) *The function $\psi \to S(\psi)$ is Fréchet differentiable at $\psi_0$ with a continuously invertible derivative $\nabla S_{\psi_0}$ (on its range).*
(d) $S(\psi_0) = 0$ *and $\hat{\psi}_n$ satisfies $S_n(\hat{\psi}_n) = o_P^*(n^{-1/2})$ and converges in outer probability to $\psi_0$.*

*Then $\sqrt{n}(\hat{\psi}_n - \psi_0) \Rightarrow -\nabla S_{\psi_0}^{-1}\mathbf{Z}$.*

In our situation, the parameter $\psi = (\boldsymbol{\theta}, \Lambda) \in \Psi = \{(\boldsymbol{\theta}, \Lambda) \colon \|\boldsymbol{\theta} - \boldsymbol{\theta}_0\| + \sup_{t\in[0,\tau]}|\Lambda(t) - \Lambda_0(t)| \leq \delta\}$ for a fixed small constant $\delta$ (note $\Psi$ is a convex set). Define a set

$$\mathcal{H} = \{(\mathbf{h}_1, h_2) \colon \|\mathbf{h}_1\| \leq 1, \|h_2\|_V \leq 1\},$$

where $\|h_2\|_V$ is the total variation of $h_2$ in $[0,\tau]$ defined as

$$\sup_{0=t_0 \leq t_1 < t_2 < \cdots < t_m = \tau} \sum_{j=1}^m |h_2(t_j) - h_2(t_{j-1})|.$$

Moreover, we let

$$S_n(\psi)(\mathbf{h}_1, h_2) = \mathbf{P}_n\{l_{\boldsymbol{\theta}}(\boldsymbol{\theta}, \Lambda)^T \mathbf{h}_1 + l_\Lambda(\boldsymbol{\theta}, \Lambda)[h_2]\},$$
$$S(\psi)(\mathbf{h}_1, h_2) = \mathbf{P}\{l_{\boldsymbol{\theta}}(\boldsymbol{\theta}, \Lambda)^T \mathbf{h}_1 + l_\Lambda(\boldsymbol{\theta}, \Lambda)[h_2]\},$$

where $l_{\boldsymbol{\theta}}(\boldsymbol{\theta}, \Lambda)$ is the first derivative of the log-likelihood function from one single subject, denoted by $l(\mathbf{O}; \boldsymbol{\theta}, \Lambda)$, with respect to $\boldsymbol{\theta}$, and $l_\Lambda(\boldsymbol{\theta}, \Lambda)[h_2]$ is the derivative of $l(\mathbf{O}; \boldsymbol{\theta}, \Lambda_\varepsilon)$ at $\varepsilon = 0$, where $\Lambda_\varepsilon(t) = \int_0^t (1 + \varepsilon h_2(s)) \, d\Lambda_0(s)$. Thus, it is easy to see that $S_n$ and $S$ are both maps from $\Psi$ to $l^\infty(\mathcal{H})$ and $\sqrt{n}\{S_n(\psi) - S(\psi)\}$ is an empirical process in the space $l^\infty(\mathcal{H})$.

According to Appendix A.2, the class

$$\mathcal{G} = \left\{ l_{\boldsymbol{\theta}}(\boldsymbol{\theta}, \Lambda)^T \mathbf{h}_1 + l_\Lambda(\boldsymbol{\theta}, \Lambda)[h_2] - l_{\boldsymbol{\theta}}(\boldsymbol{\theta}_0, \Lambda_0)^T \mathbf{h}_1 - l_\Lambda(\boldsymbol{\theta}_0, \Lambda_0)[h_2], \right.$$
$$\left. \|\boldsymbol{\theta} - \boldsymbol{\theta}_0\| + \sup_{t \in [0, \tau]} |\Lambda(t) - \Lambda_0(t)| < \delta, (\mathbf{h}_1, h_2) \in \mathcal{H} \right\}$$

is P-Donsker (cf. Section 2.1 of [24]). Moreover, Appendix A.2 also implies that

$$\sup_{(\mathbf{h}_1, h_2) \in \mathcal{H}} \mathbf{P}[l_{\boldsymbol{\theta}}(\boldsymbol{\theta}, \Lambda)^T \mathbf{h}_1 + l_\Lambda(\boldsymbol{\theta}, \Lambda)[h_2] - l_{\boldsymbol{\theta}}(\boldsymbol{\theta}_0, \Lambda_0)^T \mathbf{h}_1 - l_\Lambda(\boldsymbol{\theta}_0, \Lambda_0)[h_2]]^2 \to 0$$

when $\|\boldsymbol{\theta} - \boldsymbol{\theta}_0\| + \sup_{t \in [0, \tau]} |\Lambda(t) - \Lambda_0(t)| \to 0$. Then (a) follows from Lemma 3.3.5 of [24]. By the Donsker theorem (Section 2.5 of [24]), (b) holds as a result of Appendix A.2 and the convergence is defined in the metric space $l^\infty(\mathcal{H})$. (d) is true since $(\hat{\boldsymbol{\theta}}, \hat{\Lambda})$ maximizes $\mathbf{P}_n l(\mathbf{O}; \boldsymbol{\theta}, \Lambda)$, $(\boldsymbol{\theta}_0, \Lambda_0)$ maximizes $\mathbf{P} l(\mathbf{O}; \boldsymbol{\theta}, \Lambda)$ and $(\hat{\boldsymbol{\theta}}, \hat{\Lambda})$ converges to $(\boldsymbol{\theta}_0, \Lambda_0)$ from Theorem 3.1.

It remains to verify the conditions in (c). The proof of the first half in (c) is tedious so we defer it to Appendix A.3. We only need to prove that $\nabla S_{\psi_0}$ is continuously invertible on its range in $l^\infty(\mathcal{H})$. From Appendix A.3, $\nabla S_{\psi_0}$ can be written as follows: for any $(\boldsymbol{\theta}_1, \Lambda_1)$ and $(\boldsymbol{\theta}_2, \Lambda_2)$ in $\Psi$,

$$\begin{aligned}
&\nabla S_{\psi_0}(\boldsymbol{\theta}_1 - \boldsymbol{\theta}_2, \Lambda_1 - \Lambda_2)[\mathbf{h}_1, h_2] \\
&\qquad = (\boldsymbol{\theta}_1 - \boldsymbol{\theta}_2)^T \Omega_1[\mathbf{h}_1, h_2] + \int_0^\tau \Omega_2[\mathbf{h}_1, h_2] \, d(\Lambda_1 - \Lambda_2)(t),
\end{aligned} \quad (14)$$

where both $\Omega_1$ and $\Omega_2$ are linear operators on $\mathcal{H}$ and $\Omega = (\Omega_1, \Omega_2)$ maps $\mathcal{H} \subset R^d \times BV[0, \tau]$ to $R^d \times BV[0, \tau]$, where $BV[0, \tau]$ contains all the functions with finite total variation in $[0, \tau]$. The explicit expressions of $\Omega_1$ and $\Omega_2$ are given in Appendix A.3. From (14), we can treat $(\boldsymbol{\theta}_1 - \boldsymbol{\theta}_2, \Lambda_1 - \Lambda_2)$ as an element in $l^\infty(\mathcal{H})$ via the following definition:

$$(\boldsymbol{\theta}_1 - \boldsymbol{\theta}_2, \Lambda_1 - \Lambda_2)[\mathbf{h}_1, h_2] = (\boldsymbol{\theta}_1 - \boldsymbol{\theta}_2)^T \mathbf{h}_1 + \int_0^\tau h_2(t) \, d(\Lambda_1 - \Lambda_2)(t)$$
$$\forall (\mathbf{h}_1, h_2) \in R^d \times BV[0, \tau].$$

Then $\nabla S_{\psi_0}$ can be expanded as a linear operator from $l^\infty(\mathcal{H})$ to itself. Therefore, if we can show that there exists some positive constant $\varepsilon$ such that $\varepsilon \mathcal{H} \subset \Omega(\mathcal{H})$, then for any $(\delta\boldsymbol{\theta}, \delta\Lambda) \in l^\infty(\mathcal{H})$,

$$\|\nabla S_{\psi_0}(\delta\boldsymbol{\theta}, \delta\Lambda)\|_{l^\infty(\mathcal{H})} = \sup_{(\mathbf{h}_1, h_2) \in \mathcal{H}} \left| \delta\boldsymbol{\theta}^T \Omega_1[\mathbf{h}_1, h_2] + \int_0^\tau \Omega_2[\mathbf{h}_1, h_2]\, d\delta\Lambda(t) \right|$$
$$= \|(\delta\boldsymbol{\theta}, \delta\Lambda)\|_{l^\infty(\Omega(\mathcal{H}))} \geq \varepsilon \|(\delta\boldsymbol{\theta}, \delta\Lambda)\|_{l^\infty(\mathcal{H})}.$$

Hence, $\nabla S_{\psi_0}$ is continuously invertible.

To prove $\varepsilon \mathcal{H} \subset \Omega(\mathcal{H})$ for some $\varepsilon$ is equivalent to showing that $\Omega$ is invertible. We note from Appendix A.3 that $\Omega$ is the summation of an invertible operator and a compact operator. According to Theorem 4.25 of [20], to prove the invertibility of $\Omega$, it is sufficient to verify that $\Omega$ is one to one: if $\Omega[\mathbf{h}_1, h_2] = 0$, then by choosing $\boldsymbol{\theta}_1 - \boldsymbol{\theta}_2 = \tilde{\varepsilon} \mathbf{h}_1$ and $\Lambda_1 - \Lambda_2 = \tilde{\varepsilon} \int h_2\, d\Lambda_0$ in (14) for a small constant $\tilde{\varepsilon}$, we obtain $\nabla S_{\psi_0}(\mathbf{h}_1, \int h_2\, d\Lambda_0)[\mathbf{h}_1, h_2] = 0$. By the definition of $\nabla S_{\psi_0}$, we notice that the left-hand side is the negative information matrix in the submodel $(\boldsymbol{\theta}_0 + \varepsilon \mathbf{h}_1, \Lambda_0 + \varepsilon \int h_2\, d\Lambda_0)$. Therefore, the score function along this submodel should be zero with probability one. That is, $l_{\boldsymbol{\theta}}(\boldsymbol{\theta}_0, \Lambda_0)^T \mathbf{h}_1 + l_\Lambda(\boldsymbol{\theta}_0, \Lambda_0)[h_2] = 0$; that is, if we let $(\mathbf{h}_1^y, \mathbf{h}_1^a, \mathbf{h}_1^{\boldsymbol{\beta}}, \mathbf{h}_1^{\boldsymbol{\phi}}, \mathbf{h}_1^{\boldsymbol{\gamma}})$ be the corresponding components of $\mathbf{h}_1$ for the parameters $(\sigma_y, \text{Vec}(\boldsymbol{\Sigma}_a), \boldsymbol{\beta}, \boldsymbol{\phi}, \boldsymbol{\gamma})$, respectively, and let $\mathcal{D}_a$ be the symmetric matrix such that $\text{Vec}(\mathcal{D}_a) = \mathbf{h}_1^a$, then with probability one,

$$\begin{aligned}
0 = \int_{\mathbf{a}} G(\mathbf{a}, \mathbf{O}; \boldsymbol{\theta}_0, \Lambda_0) \\
\times \left[ \frac{\mathbf{a}^T \boldsymbol{\Sigma}_{0a}^{-1} \mathcal{D}_a \boldsymbol{\Sigma}_{0a}^{-1} \mathbf{a}}{2} - \text{Tr}(\boldsymbol{\Sigma}_{0a}^{-1} \mathcal{D}_a) - \frac{N \mathbf{h}_1^y}{\sigma_{0y}} \right. \\
+ \frac{(\mathbf{Y} - \mathbf{X}^T \boldsymbol{\beta}_0 - \widetilde{\mathbf{X}}^T \mathbf{a})^T (\mathbf{Y} - \mathbf{X}^T \boldsymbol{\beta}_0 - \widetilde{\mathbf{X}}^T \mathbf{a}) \mathbf{h}_1^y}{\sigma_{0y}^3} \\
+ \frac{\mathbf{X}(\mathbf{Y} - \mathbf{X}^T \boldsymbol{\beta}_0 - \widetilde{\mathbf{X}}^T \mathbf{a}) \mathbf{h}_1^{\boldsymbol{\beta}}}{\sigma_{0y}^2} \\
+ \Delta\{(\widetilde{\mathbf{W}}(Z) \circ \mathbf{h}_1^{\boldsymbol{\phi}})^T \mathbf{a} + \mathbf{W}(Z)^T \mathbf{h}_1^{\boldsymbol{\gamma}}\} \\
\left. - \int_0^Z e^{(\boldsymbol{\phi}_0 \circ \widetilde{\mathbf{W}}(t))^T \mathbf{a} + \mathbf{W}(t)^T \boldsymbol{\gamma}_0} \{(\widetilde{\mathbf{W}}(t) \circ \mathbf{h}_1^{\boldsymbol{\phi}})^T \mathbf{a} + \mathbf{W}(t)^T \mathbf{h}_1^{\boldsymbol{\gamma}}\}\, d\Lambda_0(t) \right] d\mathbf{a} \\
+ \int_{\mathbf{a}} G(\mathbf{a}, \mathbf{O}; \boldsymbol{\theta}_0, \Lambda_0) \left[ \Delta h_2(Z) \right. \\
\left. - \int_0^Z h_2(t) e^{(\boldsymbol{\phi}_0 \circ \widetilde{\mathbf{W}}(t))^T \mathbf{a} + \mathbf{W}(t)^T \boldsymbol{\gamma}_0}\, d\Lambda_0(t) \right] d\mathbf{a}.
\end{aligned} \tag{15}$$

Note that (15) holds with probability one, so it may not hold for any $Z \in [0, \tau]$ when $\Delta = 0$. However, if we integrate both sides from $Z$ to $\tau$ and subtract the obtained equation from (15) at $\Delta = 0$ and $Z = \tau$, it is easy to show that (15) also holds for any $Z \in [0, \tau]$ when $\Delta = 0$. Particularly, we let $\Delta = 0$ and $Z = 0$ in (15), and define

$$\mathbf{V}_a = \{\mathbf{\Sigma}_{0a}^{-1} + \widetilde{\mathbf{X}}^T \widetilde{\mathbf{X}}/\sigma_{0y}^2\}^{-1} \quad \text{and} \quad \mathbf{m}_a = \mathbf{V}_a \widetilde{\mathbf{X}}^T (\mathbf{Y} - \mathbf{X}^T \boldsymbol{\beta}_0)/\sigma_{0y}^2.$$

We obtain

$$\frac{1}{2} \text{Tr}(\mathbf{\Sigma}_{0a}^{-1} \mathcal{D}_a \mathbf{\Sigma}_{0a}^{-1} \mathbf{V}_a) + \frac{1}{2} \mathbf{m}_a^T \mathbf{\Sigma}_{0a}^{-1} \mathcal{D}_a \mathbf{\Sigma}_{0a}^{-1} \mathbf{m}_a - \frac{1}{2} \text{Tr}(\mathbf{\Sigma}_{0a}^{-1} \mathcal{D}_a) - \frac{N \mathbf{h}_1^y}{\sigma_{0y}}$$

$$+ \frac{\mathbf{h}_1^y}{\sigma_{0y}^3} \Big\{ (\mathbf{Y} - \mathbf{X}^T \boldsymbol{\beta}_0)^T (\mathbf{Y} - \mathbf{X}^T \boldsymbol{\beta}_0)$$

$$- \frac{2}{\sigma_{0y}^2} (\mathbf{Y} - \mathbf{X}^T \boldsymbol{\beta}_0)^T \widetilde{\mathbf{X}} \mathbf{V}_a \widetilde{\mathbf{X}}^T (\mathbf{Y} - \mathbf{X}^T \boldsymbol{\beta}_0)$$

$$+ \text{Tr}(\widetilde{\mathbf{X}}^T \widetilde{\mathbf{X}} \mathbf{V}_a) + \mathbf{m}_a^T \widetilde{\mathbf{X}}^T \widetilde{\mathbf{X}} \mathbf{m}_a \Big\}$$

$$+ \{(\mathbf{Y} - \mathbf{X}^T \boldsymbol{\beta}_0)^T \mathbf{X} - (\widetilde{\mathbf{X}} \mathbf{m}_a)^T \mathbf{X}\} \frac{\mathbf{h}_1^{\boldsymbol{\beta}}}{\sigma_{0y}^2} = 0.$$

Examining the coefficient for $(\mathbf{Y} - \mathbf{X}^T \boldsymbol{\beta}_0)$ gives that $\mathbf{h}_1^{\boldsymbol{\beta}} = 0$. The terms without $(\mathbf{Y} - \mathbf{X}^T \boldsymbol{\beta}_0)$ give

$$(16) \qquad -\frac{1}{2\sigma_{0y}^2} \text{Tr}(\widetilde{\mathbf{X}}^T \widetilde{\mathbf{X}} \mathbf{V}_a \mathbf{\Sigma}_{0a}^{-1} \mathcal{D}_a) - \frac{N \mathbf{h}_1^y}{\sigma_{0y}} + \frac{\mathbf{h}_1^y}{\sigma_{0y}^3} \text{Tr}(\widetilde{\mathbf{X}}^T \widetilde{\mathbf{X}} \mathbf{V}_a) = 0.$$

Moreover, the coefficients for the quadratic term $(\mathbf{Y} - \mathbf{X}^T \boldsymbol{\beta}_0)(\mathbf{Y} - \mathbf{X}^T \boldsymbol{\beta}_0)^T$ are equal to

$$(17) \qquad \begin{aligned} &\frac{1}{2\sigma_{0y}^4} \widetilde{\mathbf{X}} \mathbf{V}_a \mathbf{\Sigma}_{0a}^{-1} \mathcal{D}_a \mathbf{\Sigma}_{0a}^{-1} V_a \widetilde{\mathbf{X}}^T \\ &+ \frac{\mathbf{h}_1^y}{\sigma_{0y}^3} \left[ I - \frac{2}{\sigma_{0y}^2} \widetilde{\mathbf{X}} \mathbf{V}_a \widetilde{\mathbf{X}}^T + \frac{1}{\sigma_{0y}^4} \widetilde{\mathbf{X}} V_a \widetilde{\mathbf{X}}^T \widetilde{\mathbf{X}} \mathbf{V}_a \widetilde{\mathbf{X}}^T \right] = 0. \end{aligned}$$

Multiplying both sides of (17) by $\widetilde{\mathbf{X}}^T$ from the left and by $\widetilde{\mathbf{X}}$ from the right gives

$$\frac{1}{2\sigma_{0y}^4} \widetilde{\mathbf{X}}^T \widetilde{\mathbf{X}} \mathbf{V}_a \mathbf{\Sigma}_{0a}^{-1} \mathcal{D}_a \mathbf{\Sigma}_{0a}^{-1} V_a$$

$$+ \frac{\mathbf{h}_1^y}{\sigma_{0y}^3} \left\{ I - \frac{2}{\sigma_{0y}^2} \widetilde{\mathbf{X}}^T \widetilde{\mathbf{X}} \mathbf{V}_a + \frac{1}{\sigma_{0y}^4} \widetilde{\mathbf{X}}^T \widetilde{\mathbf{X}} \mathbf{V}_a \widetilde{\mathbf{X}}^T \widetilde{\mathbf{X}} \mathbf{V}_a \right\} = 0.$$

Since $\widetilde{\mathbf{X}}^T\widetilde{\mathbf{X}}\mathbf{V}_a/\sigma_{0y}^2 = I - \mathbf{\Sigma}_{0a}^{-1}\mathbf{V}_a$, we obtain

$$\frac{1}{2\sigma_{0y}^4}\widetilde{\mathbf{X}}^T\widetilde{\mathbf{X}}\mathbf{V}_a\mathbf{\Sigma}_{0a}^{-1}\mathcal{D}_a\mathbf{\Sigma}_{0a}^{-1}\mathbf{V}_a$$

$$+ \frac{\mathbf{h}_1^y}{\sigma_{0y}^3}\left\{I - \frac{1}{\sigma_{0y}^2}\widetilde{\mathbf{X}}^T\widetilde{\mathbf{X}}\mathbf{V}_a - \frac{1}{\sigma_{0y}^2}\widetilde{\mathbf{X}}^T\widetilde{\mathbf{X}}\mathbf{V}_a\mathbf{\Sigma}_{0a}^{-1}\mathbf{V}_a\right\} = 0.$$

Furthermore, if we multiply the above equation by $\mathbf{V}_a^{-1}\mathbf{\Sigma}_{0a}$ from the right and take the trace of the matrix, we have

$$\frac{1}{2\sigma_{0y}^4}\operatorname{Tr}(\widetilde{\mathbf{X}}^T\widetilde{\mathbf{X}}\mathbf{V}_a\mathbf{\Sigma}_{0a}^{-1}\mathcal{D}_a)$$

$$+ \frac{\mathbf{h}_1^y}{\sigma_{0y}^3}\left\{\operatorname{Tr}\left(\mathbf{V}_a^{-1}\mathbf{\Sigma}_{0a} - \frac{1}{\sigma_{0y}^2}\widetilde{\mathbf{X}}^T\widetilde{\mathbf{X}}\mathbf{\Sigma}_{0a} - \frac{1}{\sigma_{0y}^2}\widetilde{\mathbf{X}}^T\widetilde{\mathbf{X}}\mathbf{V}_a\right)\right\} = 0.$$

After substituting equation (16) into the above equation, we obtain

$$\left\{N + \operatorname{Tr}\left(\frac{1}{\sigma_{0y}^2}\widetilde{\mathbf{X}}^T\widetilde{\mathbf{X}}\mathbf{\Sigma}_{0a}\right) - \operatorname{Tr}(\mathbf{V}_a^{-1}\mathbf{\Sigma}_{0a})\right\}\mathbf{h}_1^y = 0.$$

Thus, $\mathbf{h}_1^y = 0$ based on assumption (A.7) and, moreover, from (17), $\mathcal{D}_a = 0$. Next, we let $\Delta = 0$ in (15) and obtain

$$E_\mathbf{a}\left[\exp\left\{-\int_0^Z e^{(\phi_0 \circ \widetilde{\mathbf{W}}(t))^T\mathbf{a} + \mathbf{W}(t)^T\boldsymbol{\gamma}_0}\, d\Lambda_0(t)\right\}\right.$$

$$\times \int_0^Z e^{(\phi_0 \circ \widetilde{\mathbf{W}}(t))^T\mathbf{a} + \mathbf{W}(t)^T\boldsymbol{\gamma}_0}$$

$$\left.\times \{(\widetilde{\mathbf{W}}(t) \circ \mathbf{h}_1^\phi)^T\mathbf{a} + \mathbf{W}(t)^T\mathbf{h}_1^\gamma + h_2(t)\}\, d\Lambda_0(t)\right] = 0,$$

where $\mathbf{a}$ follows a normal distribution with mean

$$\{\widetilde{\mathbf{X}}^T\widetilde{\mathbf{X}}/\sigma_{0y}^2 + \mathbf{\Sigma}_{0a}^{-1}\}\{(\mathbf{Y} - \mathbf{X}^T\boldsymbol{\beta}_0)^T\widetilde{\mathbf{X}}/\sigma_{0y}^2\}^T$$

and covariance $\{\widetilde{\mathbf{X}}^T\widetilde{\mathbf{X}}/\sigma_{0y}^2 + \mathbf{\Sigma}_{0a}^{-1}\}^{-1}$. However, for any fixed $\widetilde{\mathbf{X}}$ and $\mathbf{X}$, treating $\widetilde{\mathbf{X}}^T\mathbf{Y}$ as parameters in this normal family, $\mathbf{a}$ is a complete statistic for $\widetilde{\mathbf{X}}^T\mathbf{Y}$. Therefore,

$$\int_0^Z e^{(\phi_0 \circ \widetilde{\mathbf{W}}(t))^T\mathbf{a} + \mathbf{W}(t)^T\boldsymbol{\gamma}_0}\{(\widetilde{\mathbf{W}}(t) \circ \mathbf{h}_1^\phi)^T\mathbf{a} + \mathbf{W}(t)^T\mathbf{h}_1^\gamma + h_2(t)\}\, d\Lambda_0(t) = 0.$$

From assumption (A.9), this immediately gives $\mathbf{h}_1^\phi = 0$, $\mathbf{h}_1^\gamma = 0$ and $h_2(t) \equiv 0$.

Since conditions (a)–(d) have been proved, Theorem 3.3.1 of [24] concludes that $\sqrt{n}(\hat{\boldsymbol{\theta}} - \boldsymbol{\theta}_0, \hat{\Lambda} - \Lambda_0)$ weakly converges to a tight random element in

$l^\infty(\mathcal{H})$. Moreover, we obtain

$$
\begin{aligned}
(18) \quad & \sqrt{n}\nabla S_{\psi_0}(\hat{\boldsymbol{\theta}} - \boldsymbol{\theta}_0, \hat{\Lambda} - \Lambda_0)[\mathbf{h}_1, h_2] \\
& = \sqrt{n}(\mathbf{P}_n - \mathbf{P})\{l_{\boldsymbol{\theta}}(\boldsymbol{\theta}_0, \Lambda_0)^T \mathbf{h}_1 + l_\Lambda(\boldsymbol{\theta}_0, \Lambda_0)[h_2]\} + o_p(1),
\end{aligned}
$$

where $o_p(1)$ is a random variable which converges to zero in probability in $l^\infty(\mathcal{H})$. From (14), if we denote $(\widetilde{\mathbf{h}}_1, \widetilde{h}_2) = \Omega^{-1}(\mathbf{h}_1, h_2)$, then (18) can also be written as

$$
\begin{aligned}
(19) \quad & \sqrt{n}\bigg\{(\hat{\boldsymbol{\theta}} - \boldsymbol{\theta}_0)^T \mathbf{h}_1 + \int_0^\tau h_2(t)\, d(\hat{\Lambda} - \Lambda_0)(t)\bigg\} \\
& = \sqrt{n}(\mathbf{P}_n - \mathbf{P})\{l_{\boldsymbol{\theta}}(\boldsymbol{\theta}_0, \Lambda_0)^T \widetilde{\mathbf{h}}_1 + l_\Lambda(\boldsymbol{\theta}_0, \Lambda_0)[\widetilde{h}_2]\} + o_p(1).
\end{aligned}
$$

In other words, $\sqrt{n}(\hat{\boldsymbol{\theta}} - \boldsymbol{\theta}_0, \hat{\Lambda} - \Lambda_0)$ weakly converges to a Gaussian process in $l^\infty(\mathcal{H})$. Particularly, if we choose $h_2 = 0$ in (19), then $\hat{\boldsymbol{\theta}}^T \mathbf{h}_1$ is an asymptotic linear estimator for $\boldsymbol{\theta}_0^T \mathbf{h}_1$ with influence function being $l_{\boldsymbol{\theta}}(\boldsymbol{\theta}_0, \Lambda_0)^T \widetilde{\mathbf{h}}_1 + l_\Lambda(\boldsymbol{\theta}_0, \Lambda_0)[\widetilde{h}_2]$. Since this influence function is in the linear space spanned by the score functions for $\boldsymbol{\theta}_0$ and $\Lambda_0$, Proposition 3.3.1 in [3] concludes that the influence function is the same as the efficient influence function for $\boldsymbol{\theta}_0^T \mathbf{h}_1$; that is, $\hat{\boldsymbol{\theta}}$ is an efficient estimator for $\boldsymbol{\theta}_0$ and Theorem 3.2 has been proved.

**6. Proof of Theorem 3.3.** According to the profile likelihood theory for the semiparametric model (Theorem 1 in [18]), we need to construct an approximately least favorable submodel and verify all the conditions in that theorem. To construct the least favorable submodel, from (19) there exists a vector of functions, denoted by $\widetilde{\mathbf{h}}_2$, such that $l_{\boldsymbol{\theta}}(\boldsymbol{\theta}_0, \Lambda_0) + l_\Lambda(\boldsymbol{\theta}_0, \Lambda_0)[\widetilde{\mathbf{h}}_2]$ is the efficient score function for $\boldsymbol{\theta}_0$. Thus, the least favorable submodel at $(\boldsymbol{\theta}, \Lambda)$ is given by $\xi \mapsto (\xi, \Lambda_\xi(\boldsymbol{\theta}, \Lambda))$, where $\Lambda_\xi(\boldsymbol{\theta}, \Lambda) = \Lambda + (\xi - \boldsymbol{\theta})\int \widetilde{\mathbf{h}}_2\, d\Lambda$. For this submodel, conditions (8) and (9) in [18] hold.

We note that, in the consistency proof of Theorem 3.1, when $\hat{\boldsymbol{\theta}}$ is not necessarily the maximum likelihood estimate but $\hat{\boldsymbol{\theta}} \xrightarrow{p} \boldsymbol{\theta}_0$, the same arguments as in the proofs of (ii) and (iii) give that $\hat{\Lambda}_{\hat{\boldsymbol{\theta}}}$, which maximizes $l_n(\hat{\boldsymbol{\theta}}, \Lambda)$ over $\mathcal{Z}_n$, is bounded and its limit should be equal to $\Lambda_0$. Thus, condition (10) in [18] holds. Condition (11) in [18] can be checked straightforwardly. Furthermore, using similar arguments as in Appendix A.2, we can directly check the Donsker property of the class

$$
\left\{\nabla_\xi l(\xi, \Lambda_\xi(\boldsymbol{\theta}, \Lambda)) : \|\xi - \boldsymbol{\theta}_0\| + \|\boldsymbol{\theta} - \boldsymbol{\theta}_0\| + \sup_{t \in [0,\tau]}|\Lambda(t) - \Lambda_0(t)| < \delta\right\}
$$

and the Glivenko–Cantelli property of the class

$$
\left\{\nabla^2_{\xi\xi} l(\xi, \Lambda_\xi(\boldsymbol{\theta}, \Lambda)) : \|\xi - \boldsymbol{\theta}_0\| + \|\boldsymbol{\theta} - \boldsymbol{\theta}_0\| + \sup_{t \in [0,\tau]}|\Lambda(t) - \Lambda_0(t)| < \delta\right\}
$$

for a small constant $\delta$.

Hence, Theorem 3.3 follows from Theorem 1 and Corollaries 2 and 3 in [18].

**7. Discussion.** We have derived the asymptotic properties of the maximum likelihood estimators for joint models of repeated measurements and survival time. Our results provide a theoretical justification for the maximum likelihood estimation in such types of analysis.

Our proofs can be generalized to obtaining the asymptotic properties of the maximum likelihood estimators in other joint models, for example, the joint models of repeated measurements and multivariate survival times which were discussed by Li and Lin [17] and Huang, Zeger, Anthony and Garrett [14], the joint models of repeated measurements and recurrent event times [11], and so on. Moreover, based on our proof, we can see that Theorems 3.1–3.3 hold even if the random effect, $\mathbf{a}$, has slightly heavier tails than the normal density, for example, if the tail approximates $\exp\{-\|\mathbf{a}\|^{\alpha}\}$, where $1 < \alpha < 2$. We also note that recent work by Tsiatis and Davidian [22] used the conditional score equation approach to obtain consistent estimates for the regression coefficients in the proportional hazards model (1) without any assumptions on random effects. However, their estimates are not efficient and the applicability of the maximum likelihood estimation under this situation is yet unknown.

## APPENDIX

**A.1. Donsker property of $\mathcal{F}$.** Recall that $\mathcal{F} = \{Q(z, \mathbf{O}; \boldsymbol{\theta}, \Lambda) : z \in [0, \tau], \boldsymbol{\theta} \in \Theta, \Lambda \in \mathcal{A}\}$, where $\mathcal{A} = \{\Lambda \in \mathcal{Z}, \Lambda(\tau) \leq B_0\}$. We can rewrite $Q(z, \mathbf{O}; \boldsymbol{\theta}, \Lambda)$ as

$$Q(z, \mathbf{O}; \boldsymbol{\theta}, \Lambda) = Q_1(z, \mathbf{O}; \boldsymbol{\theta}) \frac{Q_2(z, \mathbf{O}; \boldsymbol{\theta}, \Lambda)}{Q_3(z, \mathbf{O}; \boldsymbol{\theta}, \Lambda)},$$

where

$Q_1(z, \mathbf{O}; \boldsymbol{\theta})$
$$= \exp\{\mathbf{W}(z)^T \boldsymbol{\gamma} + \Delta((\boldsymbol{\phi} \circ \widetilde{\mathbf{W}}(z)) + 2\widetilde{\mathbf{X}}^T(\mathbf{Y} - \mathbf{X}^T \boldsymbol{\beta}))^T \mathbf{V}^{-1}(\boldsymbol{\phi} \circ \widetilde{\mathbf{W}}(z))\},$$

$Q_2(z, \mathbf{O}; \boldsymbol{\theta}, \Lambda)$
$$= \int_{\mathbf{a}} \exp\left[-\frac{\mathbf{a}^T \mathbf{a}}{2} - \int_0^Z \exp\{(\boldsymbol{\phi} \circ \widetilde{\mathbf{W}}(t))^T \mathbf{V}^{-1} \mathbf{a} + \mathbf{W}(t)^T \boldsymbol{\gamma} + \mathbf{U}(t)\} d\Lambda(t)\right] d\mathbf{a},$$

$Q_3(\mathbf{O}; \boldsymbol{\theta}, \Lambda)$
$$= \int_{\mathbf{a}} \exp\left\{-\frac{\mathbf{a}^T \mathbf{a}}{2} - \int_0^Z e^{(\boldsymbol{\phi} \circ \widetilde{\mathbf{W}}(t))^T \mathbf{V}^{-1} \mathbf{a} + \mathbf{W}(t)^T \boldsymbol{\gamma}} d\Lambda(t)\right\} d\mathbf{a}.$$

Here, $\mathbf{V} = \widetilde{\mathbf{X}}^T\widetilde{\mathbf{X}}/\sigma_y^2 + \mathbf{\Sigma}_a^{-1}$ and $\mathbf{U}(t) = (\boldsymbol{\phi}\circ\widetilde{\mathbf{W}}(t))^T\mathbf{V}^{-1}((\boldsymbol{\phi}\circ\widetilde{\mathbf{W}}(z)) + \widetilde{\mathbf{X}}^T(\mathbf{Y} - \mathbf{X}^T\boldsymbol{\beta})/\sigma_y^2)$.

Using assumption (A.2), we can easily show that $Q_1(z, \mathbf{O}; \boldsymbol{\theta})$ is continuously differentiable with respect to $z$ and $\boldsymbol{\theta}$ and

$$\|\nabla_{\boldsymbol{\theta}} Q_1(z, \mathbf{O}; \boldsymbol{\theta})\| + \left|\frac{d}{dz}Q_1(z, \mathbf{O}; \boldsymbol{\theta})\right| \leq e^{g_1 + g_2\|\mathbf{Y}\|}$$

for some positive constants $g_1$ and $g_2$. Furthermore, it holds that

$$\|\nabla_{\boldsymbol{\theta}} Q_2(z, \mathbf{O}; \boldsymbol{\theta}, \Lambda)\| + \left|\frac{d}{dz}Q_2(z, \mathbf{O}; \boldsymbol{\theta}, \Lambda)\right|$$
$$\leq \int_{\mathbf{a}} \exp\left\{-\frac{\mathbf{a}^T\mathbf{a}}{2}\right\} e^{g_3\|\mathbf{a}\| + g_4\|\mathbf{Y}\| + g_5} B_0 \, d\mathbf{a}$$
$$\leq e^{g_6 + g_7\|\mathbf{Y}\|}$$

and $\|\nabla_{\boldsymbol{\theta}} Q_3(\mathbf{O}; \boldsymbol{\theta}, \Lambda)\| \leq e^{g_8 + g_9\|\mathbf{Y}\|}$ for some positive constants $g_3, \ldots, g_9$. Additionally,

$$|Q_2(z, \mathbf{O}; \boldsymbol{\theta}, \Lambda_1) - Q_2(z, \mathbf{O}; \boldsymbol{\theta}, \Lambda_2)|$$
$$\leq \left|\int_{\mathbf{a}} \exp\left\{-\frac{\mathbf{a}^T\mathbf{a}}{2}\right\} \int_0^Z e^{(\boldsymbol{\phi}\circ\widetilde{\mathbf{W}}(t))^T\mathbf{V}^{-1}\mathbf{a} + \mathbf{W}(t)^T\boldsymbol{\gamma} + \mathbf{U}(t)} \, d(\Lambda_1 - \Lambda_2)(t) \, d\mathbf{a}\right|$$
$$\leq (2\pi)^{d_a/2} \left|\int_0^Z e^{(\boldsymbol{\phi}\circ\widetilde{\mathbf{W}}(t))^T\mathbf{V}^{-2}(\boldsymbol{\phi}\circ\widetilde{\mathbf{W}}(t))/2 + \mathbf{W}(t)^T\boldsymbol{\gamma} + \mathbf{U}(t)} \, d(\Lambda_1 - \Lambda_2)(t)\right|$$
$$\leq (2\pi)^{d_a/2} \int_0^Z |\Lambda_1(t) - \Lambda_2(t)|$$
$$\qquad\times \left|\frac{d}{dt}[e^{(\boldsymbol{\phi}\circ\widetilde{\mathbf{W}}(t))^T\mathbf{V}^{-2}(\boldsymbol{\phi}\circ\widetilde{\mathbf{W}}(t))/2 + \mathbf{W}(t)^T\boldsymbol{\gamma} + \mathbf{U}(t)}]\right| dt$$
$$\qquad + (2\pi)^{d_a/2}|\Lambda_1(Z) - \Lambda_2(Z)|e^{(\boldsymbol{\phi}\circ\widetilde{\mathbf{W}}(Z))^T\mathbf{V}^{-2}(\boldsymbol{\phi}\circ\widetilde{\mathbf{W}}(Z))/2 + \mathbf{W}(Z)^T\boldsymbol{\gamma} + \mathbf{U}(Z)}$$
$$\leq e^{g_{10} + g_{11}\|\mathbf{Y}\|}\left\{|\Lambda_1(Z) - \Lambda_2(Z)| + \int_0^\tau |\Lambda_1(t) - \Lambda_2(t)| \, dt\right\},$$

where $g_{10}$ and $g_{11}$ are two positive constants. Similarly,

$$|Q_3(z, \mathbf{O}; \boldsymbol{\theta}, \Lambda_1) - Q_3(z, \mathbf{O}; \boldsymbol{\theta}, \Lambda_2)|$$
$$\leq e^{g_{10} + g_{11}\|\mathbf{Y}\|}\left\{|\Lambda_1(Z) - \Lambda_2(Z)| + \int_0^\tau |\Lambda_1(t) - \Lambda_2(t)| \, dt\right\}.$$

On the other hand, there exist positive constants $g_{13}, \ldots, g_{16}$ such that $|Q_1(z, \mathbf{O}; \boldsymbol{\theta})| \leq e^{g_{12} + g_{13}\|\mathbf{Y}\|}$, $|Q_2(z, \mathbf{O}; \boldsymbol{\theta}, \Lambda)| \leq (2\pi)^{-d_a/2}$ and $Q_3(z, \mathbf{O}; \boldsymbol{\theta}, \Lambda) \geq$

$\int_{\mathbf{a}} \exp\{-\frac{\mathbf{a}^T\mathbf{a}}{2} - e^{g_{14}+g_{15}\|\mathbf{a}\|} B_0\} \geq g_{16} > 0$. Therefore, by the mean-value theorem, we conclude that, for any $(z_1, \boldsymbol{\theta}_1, \Lambda_1)$ and $(z_2, \boldsymbol{\theta}_2, \Lambda_2)$ in $[0,\tau] \times \Theta \times \mathcal{A}$,

$$|Q(z_1, \mathbf{O}; \boldsymbol{\theta}_1, \Lambda_1) - Q(z_2, \mathbf{O}; \boldsymbol{\theta}_2, \Lambda_2)|$$
(A.1)
$$\leq e^{g_{17}+g_{18}\|\mathbf{Y}\|}\bigg\{\|\boldsymbol{\theta}_1 - \boldsymbol{\theta}_2\| + |\Lambda_1(Z) - \Lambda_2(Z)|$$
$$+ \int_0^Z |\Lambda_1(t) - \Lambda_2(t)|\,dt + |z_1 - z_2|\bigg\}$$

holds for some positive constants $g_{17}$ and $g_{18}$.

According to Theorem 2.7.5 in [24], the entropy number for the class $\mathcal{A}$ satisfies $\log N_{[\cdot]}(\varepsilon, \mathcal{A}, L_2(P)) \leq K/\varepsilon$, where $K$ is a constant. Thus, we can find $\exp\{K/\varepsilon\}$ brackets, $\{[L_j, U_j]\}$, to cover the class $\mathcal{A}$ such that, for each pair of $[L_j, U_j]$, $\|U_j - L_j\|_{L_2(P)} \leq \varepsilon$. We can further find a partition of $[0,\tau] \times \Theta$, say $S_1 \cup S_2 \cup \cdots$, such that the number of partitions is of the order $(1/\varepsilon)^{d+1}$ and for any $(z_1, \boldsymbol{\theta}_1)$ and $(z_2, \boldsymbol{\theta}_2)$ in the same partition, their Euclidean distance is less than $\varepsilon$. Therefore, the partition $\{S_1, S_2, \ldots\} \times \{[L_j, U_j]\}$ bracket covers $[0,\tau] \times \Theta \times \mathcal{A}$ and the total number of the partition is of order $(1/\varepsilon)^{d+1} \exp\{1/\varepsilon\}$. Thus, from (A.1), for any $S_k$ and $[L_j, U_j]$, the set of the functions $\{Q(z, \mathbf{O}; \boldsymbol{\theta}, \Lambda) : (z, \boldsymbol{\theta}) \in S_k, \Lambda \in \mathcal{A}, \Lambda \in [L_j, U_j]\}$ can be bracket covered by

$$\bigg[Q(z_k, \mathbf{O}; \boldsymbol{\theta}_k, \Lambda_j)$$
$$- e^{g_{17}+g_{18}\|\mathbf{Y}\|}\bigg\{\varepsilon + |U_j(Z) - L_j(Z)| + \int_0^\tau |U_j(t) - L_j(t)|\,dt\bigg\},$$
$$Q(z_k, \mathbf{O}; \boldsymbol{\theta}_k, \Lambda_j)$$
$$+ e^{g_{17}+g_{18}\|\mathbf{Y}\|}\bigg\{\varepsilon + |U_j(Z) - L_j(Z)| + \int_0^\tau |U_j(t) - L_j(t)|\,dt\bigg\}\bigg],$$

where $(z_k, \boldsymbol{\theta}_k)$ is a fixed point in $S_k$ and $\Lambda_j$ is a fixed function in $[L_j, U_j]$. Note that the $L_2(P)$ distance between these two functions is less than $O(\varepsilon)$. Therefore, we have

$$N_{[\cdot]}(\varepsilon, \mathcal{F}, \|\cdot\|_{L_2(P)}) \leq O(1)\left(\frac{1}{\varepsilon}\right)^{d+1} e^{1/\varepsilon}.$$

Furthermore, $\mathcal{F}$ has an $L_2(P)$-integrable covering function, which is equal to $O(e^{g_{17}+g_{18}\|\mathbf{Y}\|})$. From Theorem 2.5.6 in [24], $\mathcal{F}$ is P-Donsker.

In the above derivation, we also note that all the functions in $\mathcal{F}$ are bounded from below by $e^{-g_{19}-g_{20}\|\mathbf{Y}\|}$ for some positive constants $g_{19}$ and $g_{20}$.

**A.2. Donsker property of $\mathcal{G}$.** For any $(\mathbf{h}_1, h_2) \in \mathcal{H}$, if we let $(\mathbf{h}_1^y, \mathbf{h}_1^a, \mathbf{h}_1^{\boldsymbol{\beta}}, \mathbf{h}_1^{\boldsymbol{\phi}}, \mathbf{h}_1^{\boldsymbol{\gamma}})$ be the corresponding components of $\mathbf{h}_1$ for the parameters $(\sigma_y, \boldsymbol{\Sigma}_a, \boldsymbol{\beta}, \boldsymbol{\phi}, \boldsymbol{\gamma})$, respectively, then $l_{\boldsymbol{\theta}}(\boldsymbol{\theta}, \Lambda)^T \mathbf{h}_1 + l_\Lambda(\boldsymbol{\theta}, \Lambda)[h_2]$ has the expression

$$\left[\mu_1(\mathbf{O}; \boldsymbol{\theta}, \Lambda)^T \mathbf{h}_1 - \int_0^Z \mu_2(t, \mathbf{O}; \boldsymbol{\theta}, \Lambda)^T \mathbf{h}_1 \, d\Lambda(t)\right]$$

$$+ \Delta h_2(Z) - \int_0^Z \mu_3(t, \mathbf{O}; \boldsymbol{\theta}, \Lambda) h_2(t) \, d\Lambda(t),$$

where

$$\mu_1(\mathbf{O}; \boldsymbol{\theta}, \Lambda)^T \mathbf{h}_1$$

$$= \left\{\int_{\mathbf{a}} G(\mathbf{a}, \mathbf{O}; \boldsymbol{\theta}, \Lambda) \, d\mathbf{a}\right\}^{-1}$$

$$\times \int_{\mathbf{a}} G(\mathbf{a}, \mathbf{O}; \boldsymbol{\theta}, \Lambda) \left[\frac{\mathbf{a}^T \boldsymbol{\Sigma}_a^{-1} \mathcal{D}_a \boldsymbol{\Sigma}_a^{-1} \mathbf{a}}{2} - \text{Tr}(\boldsymbol{\Sigma}_a^{-1} \mathcal{D}_a) - \frac{N \mathbf{h}_1^y}{\sigma_y}\right.$$

$$+ \frac{(\mathbf{Y} - \mathbf{X}^T \boldsymbol{\beta} - \widetilde{\mathbf{X}}^T \mathbf{a})^T (\mathbf{Y} - \mathbf{X}^T \boldsymbol{\beta} - \widetilde{\mathbf{X}}^T \mathbf{a}) \mathbf{h}_1^y}{\sigma_y^3}$$

$$+ \frac{\mathbf{X}(\mathbf{Y} - \mathbf{X}^T \boldsymbol{\beta} - \widetilde{\mathbf{X}}^T \mathbf{a}) \mathbf{h}_1^{\boldsymbol{\beta}}}{\sigma_y^2}$$

$$\left. + \Delta\{(\widetilde{\mathbf{W}}(Z) \circ \mathbf{h}_1^{\boldsymbol{\phi}})^T \mathbf{a} + \mathbf{W}(Z)^T \mathbf{h}_1^{\boldsymbol{\gamma}}\}\right] d\mathbf{a},$$

$$\mu_2(t, \mathbf{O}; \theta, \Lambda)^T \mathbf{h}_1$$

$$= \left\{\int_{\mathbf{a}} G(\mathbf{a}, \mathbf{O}; \boldsymbol{\theta}, \Lambda) \, d\mathbf{a}\right\}^{-1}$$

$$\times \int_{\mathbf{a}} G(\mathbf{a}, \mathbf{O}; \boldsymbol{\theta}, \Lambda) e^{(\boldsymbol{\phi} \circ \widetilde{\mathbf{W}}(t))^T \mathbf{a} + \mathbf{W}(t)^T \boldsymbol{\gamma}}$$

$$\times \{(\widetilde{\mathbf{W}}(t) \circ \mathbf{h}_1^{\boldsymbol{\phi}})^T \mathbf{a} + \mathbf{W}(t)^T \mathbf{h}_1^{\boldsymbol{\gamma}}\} \, d\mathbf{a}$$

and

$$\mu_3(t, \mathbf{O}; \boldsymbol{\theta}, \Lambda)$$

$$= \left\{\int_{\mathbf{a}} G(\mathbf{a}, \mathbf{O}; \boldsymbol{\theta}, \Lambda) \, d\mathbf{a}\right\}^{-1} \int_{\mathbf{a}} G(\mathbf{a}, \mathbf{O}; \boldsymbol{\theta}, \Lambda) e^{(\boldsymbol{\phi} \circ \widetilde{\mathbf{W}}(t))^T \mathbf{a} + \mathbf{W}(t)^T \boldsymbol{\gamma}} \, d\mathbf{a}.$$

Here, $\mathcal{D}_a$ is a symmetric matrix such that $\text{Vec}(\mathcal{D}_a) = \mathbf{h}_1^a$.

For $j = 1, 2, 3$, we denote $\nabla_{\boldsymbol{\theta}} \mu_j$ and $\nabla_\Lambda \mu_j[\delta\Lambda]$ as the derivatives of $\mu_j$ with respect to $\boldsymbol{\theta}$ and $\Lambda$ along the path $\Lambda + \varepsilon \delta\Lambda$. Then using similar calculations to Appendix A.1, it is tedious to verify that $\nabla_\Lambda \mu_j[\delta\Lambda] = \int_0^t \mu_{j+3}(s, \mathbf{O}; \boldsymbol{\theta}, \Lambda) \, d\delta\Lambda(s)$ and that there exist two positive constants $r_1$ and $r_2$ such that

$$\sum_j \{|\mu_j| + |\nabla_{\boldsymbol{\theta}} \mu_j|\} \leq e^{r_1 + r_2 \|b\mathbf{Y}\|}.$$

On the other hand, by the mean value theorem, we have that, for any $(\boldsymbol{\theta}, \Lambda, \mathbf{h}_1, h_2)$ and $(\widetilde{\boldsymbol{\theta}}, \widetilde{\Lambda}, \widetilde{\mathbf{h}}_1, \widetilde{h}_2)$ in $\Psi \times \mathcal{H}$,

$$l_{\boldsymbol{\theta}}(\boldsymbol{\theta}, \Lambda)^T \mathbf{h}_1 + l_\Lambda(\boldsymbol{\theta}, \Lambda)[h_2] - l_{\boldsymbol{\theta}}(\widetilde{\boldsymbol{\theta}}, \widetilde{\Lambda})^T \widetilde{\mathbf{h}}_1 - l_\Lambda(\widetilde{\boldsymbol{\theta}}, \widetilde{\Lambda})[\widetilde{h}_2]$$

$$= (\boldsymbol{\theta} - \widetilde{\boldsymbol{\theta}})^T \nabla_{\boldsymbol{\theta}} \mu_1(\mathbf{O}; \boldsymbol{\theta}^*, \Lambda^*) \mathbf{h}_1 + \int_0^Z \mu_4(t, \mathbf{O}; \boldsymbol{\theta}^*, \Lambda^*)^T \mathbf{h}_1 \, d(\Lambda - \widetilde{\Lambda})(t)$$

$$- \int_0^Z (\boldsymbol{\theta} - \widetilde{\boldsymbol{\theta}})^T \nabla_{\boldsymbol{\theta}} \mu_2(t, \mathbf{O}; \boldsymbol{\theta}^*, \Lambda^*) \mathbf{h}_1 \, d\Lambda(t)$$

$$- \int_0^Z \int_0^t \mu_5(s, \mathbf{O}; \boldsymbol{\theta}^*, \Lambda^*) \, d(\Lambda - \widetilde{\Lambda})(s) \mathbf{h}_1 \, d\Lambda(t)$$

$$- \int_0^Z \mu_2(t, \mathbf{O}; \boldsymbol{\theta}^*, \Lambda^*)^T \mathbf{h}_1 \, d(\Lambda - \widetilde{\Lambda})(t)$$

(A.2)
$$- \int_0^Z (\boldsymbol{\theta} - \widetilde{\boldsymbol{\theta}})^T \nabla_{\boldsymbol{\theta}} \mu_3(t, \mathbf{O}; \boldsymbol{\theta}^*, \Lambda^*) h_2(t) \, d\Lambda(t)$$

$$- \int_0^Z \int_0^t \mu_6(s, \mathbf{O}; \boldsymbol{\theta}^*, \Lambda^*) \, d(\Lambda - \widetilde{\Lambda})(s) h_2(t) \, d\Lambda(t)$$

$$- \int_0^Z \mu_3(t, \mathbf{O}; \boldsymbol{\theta}^*, \Lambda^*) h_2(t) \, d(\Lambda - \widetilde{\Lambda})(t)$$

$$+ \mu_1(\mathbf{O}; \widetilde{\boldsymbol{\theta}}, \widetilde{\Lambda})^T (\mathbf{h}_1 - \widetilde{\mathbf{h}}_1) - \int_0^Z \mu_2(t, \mathbf{O}; \widetilde{\boldsymbol{\theta}}, \widetilde{\Lambda})^T (\mathbf{h}_1 - \widetilde{\mathbf{h}}_1) \, d\Lambda(t)$$

$$+ \Delta(h_2(Z) - \widetilde{h}_2(Z)) - \int_0^Z \mu_3(t, \mathbf{O}; \widetilde{\boldsymbol{\theta}}, \widetilde{\Lambda})(h_2(t) - \widetilde{h}_2(t)) \, d\widetilde{\Lambda}(t),$$

where $(\boldsymbol{\theta}^*, \Lambda^*)$ is equal to $\varepsilon^*(\boldsymbol{\theta}, \Lambda) + (1 - \varepsilon^*)(\widetilde{\boldsymbol{\theta}}, \widetilde{\Lambda})$ for some $\varepsilon^* \in [0, 1]$. Thus,

$$|l_{\boldsymbol{\theta}}(\boldsymbol{\theta}, \Lambda)^T \mathbf{h}_1 + l_\Lambda(\boldsymbol{\theta}, \Lambda)[h_2] - l_{\boldsymbol{\theta}}(\widetilde{\boldsymbol{\theta}}, \widetilde{\Lambda})^T \widetilde{\mathbf{h}}_1 - l_\Lambda(\widetilde{\boldsymbol{\theta}}, \widetilde{\Lambda})[\widetilde{h}_2]|$$

$$\leq e^{r_1 + r_2 \|\mathbf{Y}\|} \bigg\{ \|\boldsymbol{\theta} - \widetilde{\boldsymbol{\theta}}\| + \|\mathbf{h}_1 - \widetilde{\mathbf{h}}_1\| + |\Lambda(Z) - \widetilde{\Lambda}(Z)|$$

$$+ \int_0^\tau |\Lambda(t) - \widetilde{\Lambda}(t)|[dt + d|h_2(t)| + d|\widetilde{h}_2(t)|]$$

$$+ |h_2(Z) - \widetilde{h}_2(Z)| + \int_0^\tau |h_2(t) - \widetilde{h}_2(t)|[d\Lambda_1(t) + d\Lambda_2(t)] \bigg\},$$

where $d|h_2(t)| = dh_2^+(t) + dh_2^-(t)$ and $d|\widetilde{h}_2(t)| = d\widetilde{h}_2^+(t) + d\widetilde{h}_2^-(t)$. Therefore, by using the same arguments as in Appendix A.1 and noting that $\log N_{[\,\cdot\,]}(\varepsilon, \{h_2 : \|h_2\|_V \leq B_1\}, L_2(Q)) \leq K/\varepsilon$ for a constant $B_1$ and any probability measure $Q$ where $K$ is a constant (Theorem 2.7.5 of [24]), we obtain

$$\log N_{[\,\cdot\,]}(\varepsilon, \mathcal{G}, L_2(P)) \leq O\Big(\frac{1}{\varepsilon} + \log \varepsilon\Big).$$

Hence, $\mathcal{G}$ is P-Donsker.

Furthermore, from (A.2) we can calculate that

$$|l_{\boldsymbol{\theta}}(\boldsymbol{\theta}, \Lambda)^T \mathbf{h}_1 + l_\Lambda(\boldsymbol{\theta}, \Lambda)[h_2] - l_{\boldsymbol{\theta}}(\boldsymbol{\theta}_0, \Lambda_0)^T \mathbf{h}_1 - l_\Lambda(\boldsymbol{\theta}_0, \Lambda_0)[h_2]|$$
$$\leq e^{r_1 + r_2 \|\mathbf{Y}\|} \Big\{ \|\boldsymbol{\theta} - \boldsymbol{\theta}_0\| + |\Lambda(Z) - \Lambda_0(Z)| + \int_0^\tau |\Lambda(t) - \Lambda_0(t)| \, dt \Big\}$$
$$+ \Big| \int_0^Z \mu_3(t, \mathbf{O}; \boldsymbol{\theta}^*, \Lambda^*) h_2(t) \, d(\Lambda - \Lambda_0)(t) \Big|.$$

If $\|\boldsymbol{\theta} - \boldsymbol{\theta}_0\| \to 0$ and $\sup_{t \in [0, \tau]} |\Lambda(t) - \Lambda_0(t)| \to 0$, the above expression converges to zero uniformly. Thus,

$$\sup_{(\mathbf{h}_1, h_2) \in \mathcal{H}} \mathbf{P}[l_{\boldsymbol{\theta}}(\boldsymbol{\theta}, \Lambda)^T \mathbf{h}_1 + l_\Lambda(\boldsymbol{\theta}, \Lambda)[h_2] - l_{\boldsymbol{\theta}}(\boldsymbol{\theta}_0, \Lambda_0)^T \mathbf{h}_1 - l_\Lambda(\boldsymbol{\theta}_0, \Lambda_0)[h_2]]^2 \to 0.$$

**A.3. Derivative operator $\nabla S_{\psi_0}$.** From (A.2) we can obtain that

$$l_{\boldsymbol{\theta}}(\boldsymbol{\theta}, \Lambda)^T \mathbf{h}_1 + l_\Lambda(\boldsymbol{\theta}, \Lambda)[h_2] - l_{\boldsymbol{\theta}}(\boldsymbol{\theta}_0, \Lambda_0)^T \mathbf{h}_1 - l_\Lambda(\boldsymbol{\theta}_0, \Lambda_0)[h_2]$$
$$= (\boldsymbol{\theta} - \boldsymbol{\theta}_0)^T \Big\{ \nabla_{\boldsymbol{\theta}} \mu_1(\mathbf{O}; \boldsymbol{\theta}^*, \Lambda^*) - \int_0^Z \nabla_{\boldsymbol{\theta}} \mu_2(t, \mathbf{O}; \boldsymbol{\theta}^*, \Lambda^*) \, d\Lambda_0(t) \Big\} \mathbf{h}_1$$
$$+ \mathbf{h}_1^T \int_0^\tau I(t \leq Z) \Big\{ \mu_4(t, \mathbf{O}; \boldsymbol{\theta}^*, \Lambda^*) - \mu_2(t, \mathbf{O}; \boldsymbol{\theta}^*, \Lambda^*)$$
$$- \mu_5(t, \mathbf{O}; \boldsymbol{\theta}^*, \Lambda^*) \int_t^Z d\Lambda_0(s) \Big\} d(\Lambda - \Lambda_0)(t)$$
$$- (\boldsymbol{\theta} - \boldsymbol{\theta}_0)^T \int_0^\tau I(t \leq Z) \nabla_{\boldsymbol{\theta}} \mu_3(t, \mathbf{O}; \boldsymbol{\theta}^*, \Lambda^*) h_2(t) \, d\Lambda_0(t)$$
$$- \int_0^\tau \Big\{ I(t \leq Z) \mu_6(t, \mathbf{O}; \boldsymbol{\theta}^*, \Lambda^*) \int_t^Z h_2(s) \, d\Lambda_0(s)$$
$$+ I(t \leq Z) \mu_3(t, \mathbf{O}; \boldsymbol{\theta}^*, \Lambda^*) h_2(t) \Big\} d(\Lambda - \Lambda_0)(t).$$

Then it is clear that

$$\nabla S_{\psi_0}(\boldsymbol{\theta} - \boldsymbol{\theta}_0, \Lambda - \Lambda_0)[\mathbf{h}_1, h_2]$$
$$= (\boldsymbol{\theta} - \boldsymbol{\theta}_0)^T \mathbf{P} \Big\{ \nabla_{\boldsymbol{\theta}} \mu_1(\mathbf{O}; \boldsymbol{\theta}_0, \Lambda_0) - \int_0^Z \nabla_{\boldsymbol{\theta}} \mu_2(t, \mathbf{O}; \boldsymbol{\theta}_0, \Lambda_0) \, d\Lambda_0(t) \Big\} \mathbf{h}_1$$

$$+ \mathbf{h}_1^T \int_0^\tau \mathbf{P}\bigg[I(t \leq Z)\bigg\{\mu_4(t,\mathbf{O};\boldsymbol{\theta}_0,\Lambda_0)$$
$$- \mu_2(t,\mathbf{O};\boldsymbol{\theta}_0,\Lambda_0)$$
$$- \mu_5(t,\mathbf{O};\boldsymbol{\theta}_0,\Lambda_0)\int_t^Z d\Lambda_0(s)\bigg\}\bigg] d(\Lambda - \Lambda_0)(t)$$
$$- (\boldsymbol{\theta} - \boldsymbol{\theta}_0)^T \int_0^\tau \mathbf{P}\{I(t \leq Z)\nabla_{\boldsymbol{\theta}}\mu_3(t,\mathbf{O};\boldsymbol{\theta}_0,\Lambda_0)\}h_2(t)\,d\Lambda_0(t)$$
$$- \int_0^\tau \mathbf{P}\bigg\{I(t \leq Z)\mu_6(t,\mathbf{O};\boldsymbol{\theta}_0,\Lambda_0)\int_t^Z h_2(s)\,d\Lambda_0(s)$$
$$+ I(t \leq Z)\mu_3(t,\mathbf{O};\boldsymbol{\theta}_0,\Lambda_0)h_2(t)\bigg\} d(\Lambda - \Lambda_0)(t).$$

It is tedious to check that, for $j = 1, 2, \ldots, 6$,
$$\sup_{t \in [0,\tau]} \|\mu_j(t,\mathbf{O};\boldsymbol{\theta}^*,\Lambda^*) - \mu_j(t,\mathbf{O};\boldsymbol{\theta}_0,\Lambda_0)\|$$
$$\leq e^{r_3 + r_4\|\mathbf{Y}\|}\bigg\{\|\boldsymbol{\theta}^* - \boldsymbol{\theta}_0\| + \sup_{t \in [0,\tau]}|\Lambda^*(t) - \Lambda_0(t)|\bigg\}.$$

Thus,
$$\mathbf{P}[l_{\boldsymbol{\theta}}(\boldsymbol{\theta},\Lambda)^T\mathbf{h}_1 + l_\Lambda(\boldsymbol{\theta},\Lambda)[h_2] - l_{\boldsymbol{\theta}}(\boldsymbol{\theta}_0,\Lambda_0)^T\mathbf{h}_1 - l_\Lambda(\boldsymbol{\theta}_0,\Lambda_0)[h_2]]$$
$$= \nabla S_{\psi_0}(\boldsymbol{\theta} - \boldsymbol{\theta}_0, \Lambda - \Lambda_0)[\mathbf{h}_1, h_2]$$
$$+ o\bigg(\|\boldsymbol{\theta} - \boldsymbol{\theta}_0\| + \sup_{t \in [0,\tau]}|\Lambda(t) - \Lambda_0(t)|\bigg)(\|\mathbf{h}_1\| + \|h_2\|_V).$$

Therefore, $S(\psi_0)$ is Fréchet differentiable.

We can rewrite $\nabla S_{\psi_0}(\boldsymbol{\theta} - \boldsymbol{\theta}_0, \Lambda - \Lambda_0)[\mathbf{h}_1, h_2]$ as $\mathbf{h}_1^T\Omega_1[\mathbf{h}_1, h_2] + \int_0^\tau \Omega_2[\mathbf{h}_1, h_2]\,d(\Lambda - \Lambda_0)$, where
$$\Omega_1[\mathbf{h}_1, h_2] = \mathbf{h}_1^T \mathbf{P}\bigg\{\nabla_{\boldsymbol{\theta}}\mu_1(\mathbf{O};\boldsymbol{\theta}_0,\Lambda_0) - \int_0^Z \nabla_{\boldsymbol{\theta}}\mu_2(t,\mathbf{O};\boldsymbol{\theta}_0,\Lambda_0)\,d\Lambda_0(t)\bigg\}$$
$$- \int_0^\tau \mathbf{P}\{I(t \leq Z)\nabla_{\boldsymbol{\theta}}\mu_3(t,\mathbf{O};\boldsymbol{\theta}_0,\Lambda_0)\}h_2(t)\,d\Lambda_0(t)$$

and
$$\Omega_2[\mathbf{h}_1, h_2] = \mathbf{h}_1^T \mathbf{P}\bigg[I(t \leq Z)\bigg\{\mu_4(t,\mathbf{O};\boldsymbol{\theta}_0,\Lambda_0) - \mu_2(t,\mathbf{O};\boldsymbol{\theta}_0,\Lambda_0)$$
$$- \mu_5(t,\mathbf{O};\boldsymbol{\theta}_0,\Lambda_0)\int_t^Z d\Lambda_0(s)\bigg\}\bigg]$$
$$- \mathbf{P}\bigg\{I(t \leq Z)\mu_6(t,\mathbf{O};\boldsymbol{\theta}_0,\Lambda_0)\int_t^Z h_2(s)\,d\Lambda_0(s)\bigg\}$$

$$- \mathbf{P}\{I(t \leq Z)\mu_3(t, \mathbf{O}; \boldsymbol{\theta}_0, \Lambda_0)\}h_2(t).$$

Then the operator $\Omega = (\Omega_1, \Omega_2)$ is the bounded linear operator from $R^d \times BV[0, \tau]$ to itself. Moreover, we note that $\Omega = \mathbf{A} + (\mathbf{K}_1, \mathbf{K}_2)$, where $\mathbf{A}(\mathbf{h}_1, h_2) = (\mathbf{h}_1, -\mathbf{P}\{I(t \leq Z)\mu_3(t, \mathbf{O}; \boldsymbol{\theta}_0, \Lambda_0)\}h_2(t))$, $\mathbf{K}_1(\mathbf{h}_1, h_2) = \Omega_1[\mathbf{h}_1, h_2] - \mathbf{h}_1$ and

$$\mathbf{K}_2(\mathbf{h}_1, h_2) = \mathbf{h}_1^T \mathbf{P}\bigg[ I(t \leq Z)\bigg\{ \mu_4(t, \mathbf{O}; \boldsymbol{\theta}_0, \Lambda_0) - \mu_2(t, \mathbf{O}; \boldsymbol{\theta}_0, \Lambda_0) - \mu_5(t, \mathbf{O}; \boldsymbol{\theta}_0, \Lambda_0) \int_t^Z d\Lambda_0(s) \bigg\} \bigg] - \mathbf{P}\bigg\{ I(t \leq Z)\mu_6(t, \mathbf{O}; \boldsymbol{\theta}_0, \Lambda_0) \int_t^Z h_2(s) \, d\Lambda_0(s) \bigg\}.$$

Obviously, $\mathbf{A}$ is invertible. $\mathbf{K}_1$ maps into a finite-dimensional space, so it is compact. The image of $\mathbf{K}_2$ is a continuously differentiable function in $[0, \tau]$. According to the Arzela–Ascoli theorem (page 245 in [21]), $\mathbf{K}_2$ is a compact operator from $R^d \times BV[0, \tau]$ to $BV[0, \tau]$. Therefore, we conclude that $\Omega$ is the summation of an invertible operator and a compact operator.

**Acknowledgments.** We thank the anonymous reviewers for their valuable suggestions, which lead to great improvement from our original version. Our thanks also goes to Dr. Clarence E. Davis, who has carefully read our manuscript and helped improve the presentation of this paper.

Department of Biostatistics
University of North Carolina
Chapel Hill, North Carolina 27599-7420
USA
e-mail: dzeng@bios.unc.edu
cai@bios.unc.edu